\documentclass{article}

\begin{document}

\title{Deformations of multivalued harmonic functions}
\author{Simon Donaldson}
\maketitle
\newcommand{\sthird}{{\textstyle \frac{1}{3}}}
\newcommand{\shalf}{{\textstyle \frac{1}{2}}}
\newcommand{\sthreetwo}{{\textstyle \frac{3}{2}}}
\newcommand{\stwothirds}{{\textstyle \frac{2}{3}}}
\newcommand{\squart}{{\textstyle \frac{1}{4}}}
\newcommand{\sfifth}{{\textstyle \frac{1}{5}}}
\newcommand{\sfivetwo}{{\textstyle \frac{5}{2}}}
\newcommand{\sthreequart}{{\textstyle \frac{3}{4}}}
\newcommand{\tQ}{\tilde{Q}}
\newcommand{\cF}{{\cal F}}
\newcommand{\cZ}{{\cal Z}}
\newcommand{\cU}{{\cal U}}
\newcommand{\bR}{{\bf R}}
\newcommand{\bC}{{\bf C}}
\newcommand{\real}{{\rm Re}}
\newcommand{\cC}{{\cal C}}
\newcommand{\cS}{{\cal S}}
\newcommand{\cM}{{\cal M}}
\newcommand{\cP}{{\cal P}}
\newcommand{\cG}{{\cal G}}
\newcommand{\cT}{{\cal T}}
\newcommand{\tDelta}{\tilde{\Delta}}
\newcommand{\ddt}{\frac{d}{dt}}
\newcommand{\shar}{\tilde{\phi}}
\newcommand{\Lie}{{\rm Lie}}
\newcommand{\thalf}{1/2}
\newcommand{\cD}{{\cal D}}
\newcommand{\tG}{\tilde{G}}
\newcommand{\KN}{KN}
\newcommand{\tM}{\widetilde{M}}
\newcommand{\uG}{\underline{G}}
\newcommand{\cL}{{\cal L}}
\newcommand{\cH}{{\cal H}}
\newcommand{\loc}{{\rm loc}}
\newcommand{\scaled}{{\rm scaled}}
\newcommand{\uZ}{\underline{Z}}
\newcommand{\uE}{\underline{E}}
\newcommand{\Diff}{{\rm Diff}}
\newcommand{\tH}{\tilde{H}}
\newcommand{\uH}{\underline{H}}
\newcommand{\ucT}{\underline{{\cal T}}}
\newcommand{\utG}{\underline{\tilde{G}}}
\newcommand{\uL}{\underline{L}}
\newcommand{\dive}{{\rm div}}
\newcommand{\uA}{\underline{A}}
\newcommand{\bA}{{\bf A}}
\newcommand{\bB}{{\bf B}}
\newcommand{\bone}{{\bf 1}}
\newcommand{\cA}{{\cal A}}
\newtheorem{lem}{Lemma}
\newtheorem{thm}{Theorem}
\newtheorem{prop}{Proposition}
\newtheorem{cor}{Corollary}



\section{Introduction}

Multivalued, or branched, holomorphic functions are familiar in complex analysis and Riemann surface theory and  the real parts of these yield multivalued harmonic functions. In this paper we extend some of these well-known ideas to general Riemannian manifolds, studying harmonic functions which are locally modelled on $\real(z^{\shalf})$ where $z$ is a local complex co-ordinate transverse to a codimension-2 branch set.

To set up our problem precisely, let $(M,g)$ be a compact Riemannian $n$-manifold and  $\Sigma\subset M$  a codimension-2 submanifold. For simplicity we assume that $\Sigma$ is co-oriented in $M$. Write $\Gamma$ for the group of isometries of the real line, so there is an exact sequence
$$   0\rightarrow (\bR,+)\rightarrow \Gamma\rightarrow \{{\pm 1}\}\rightarrow 1 .$$

 We suppose that we have a representation $\chi:\pi_{1}(M\setminus \Sigma) \rightarrow \Gamma$. This representation  defines a flat bundle $E_{+}$ over $M\setminus \Sigma$ with fibre $\bR$ and structure group $\Gamma$. The composite of $\chi$ with the homomorphism to $\{\pm 1\}$ defines a flat vector bundle $E=E_{\chi}$,  which is the vertical tangent bundle of $E_{+}$. We assume that the  flat bundles $E, E_{+}$ are not isomorphic, {\it i.e.} that $E_{+}$ does not have a parallel global section. We also assume that $\chi$ maps any small loop $\delta$ locally linking $\Sigma$ to a reflection (i.e an element of order $2$ in $\Gamma$). This implies that if $U\subset M$ is a tubular neighbourhood of $\Sigma$ the restriction of $E_{+}$ to $U\setminus \Sigma$ is naturally a vector bundle and sections of $E_{+}$ over this neighbourhood can be viewed as 2-valued functions, which change sign as we move around the loop $\delta$ in the familiar fashion.  

Given this set-up we consider  {\it harmonic} sections $\phi$ of $E_{+}$ over $M\setminus \Sigma$, satisfying the Laplace equation $\Delta_{g}\phi=0$ which makes sense in an obvious way due to the flat structure. (More precisely, for any section $\phi$ of $E_{+}$ the Laplacian $\Delta_{g}\phi$ is a section of $E$.)  The first basic fact we need is that there is a  unique  harmonic section  $\phi$ with derivative in $L^{2}$ (see Section 4 below). The second basic fact we need, which we explain in detail in Sections 2 and 3 below,  is that this section $\phi$ has an asymptotic description near $\Sigma$. Working first near  a fixed point $p$ in $\Sigma$ and with a suitable complex co-ordinate $z$ on a slice through $p$ transverse to $\Sigma$,  the leading term  has the form
$$  \phi= \real( a z^{\shalf}) + O(\vert z\vert^{\sthreetwo}), $$
for a complex number $a$. If this leading term vanishes we have
$$  \phi= \real( b z^{\sthreetwo}) + O(\vert z\vert^{\sfivetwo}), $$ for another complex number $b$. (One can define the second term $b$ even when $a$ does not vanish  but if the mean curvature of $\Sigma$ is not zero there is a small subtlety involved, which we discuss in Section 3 below.) For a global version of  this,  let $N$ be the normal bundle of $\Sigma$ in $M$, regarded as a complex line bundle using the co-orientation. The representation $\chi$ defines a square root $N^{1/2}$, with dual $N^{-1/2}$. The global version of the leading term in the asymptotic expansion  is a section $\bA$  of $N^{-1/2}$ over $\Sigma$. If this vanishes the next term is a section  $\bB$ of $N^{-3/2}$. In sum, the harmonic section $\phi$ defines  $\bA= \bA(\Sigma,g,\chi)\in \Gamma(N^{-1/2})$ and if $\bA=0$ we have $\bB= \bB(\Sigma ,g,\chi)\in \Gamma(N^{-3/2})$.

The problem we consider in this paper is that of choosing $\Sigma$ so that the leading term $\bA(\Sigma,g,\chi)$ vanishes. This problem has a similar character to a free boundary value problem, with the difference that $\Sigma$ has codimension 2 rather than 1.  
We prove a deformation result in this direction. Let ${\cal S}$ be the space of codimension-2 submanifolds in $M$ and ${\cal M}$  be the space of Riemannian metrics on $M$. These have $C^{\infty}$ topologies in a standard way. Let $R$ be the space of representations $\rho$---as we will see below,  each connected component of $R$ can be identified with a real vector space.  
\begin{thm}
Suppose that $\bA(\Sigma_{0}, g_{0},\chi_{0})=0$ and that $\bB(\Sigma_{0}, g_{0}, \chi_{0})$ is nowhere-vanishing on $\Sigma_{0}$. Then there is a neighbourhood ${\cal U}$ of $(g_{0}, \chi_{0})$ in ${\cal M}\times R $ and a neighbourhood ${\cal V}$ of $\Sigma_{0}$ in $\cS$ such that for any $(g,\chi)$ in ${\cal U}$ there is a unique $\Sigma$ in ${\cal V}$ with $\bA(\Sigma,g,\chi)=0$.
\end{thm}
In short, the equation $\bA(\Sigma, g,\chi)=0$ locally defines $\Sigma$ implicitly in terms of $(g,\chi)$,  provided that $\bB$ is nowhere vanishing. (Strictly speaking, the space of representations $R$ depends on $\Sigma$ but since we are only considering small deformations these can be identified, in an obvious way.)

\

We now give five items of background and motivation for this study. 
\begin{enumerate}
\item The data $(M,\Sigma,\chi)$ has an alternative description in terms of a branched cover. The composite of $\chi$ with $\Gamma\rightarrow \{\pm 1\}$ defines a double cover of $M\setminus \Sigma$ and our assumptions  imply that this extends to a branched cover
$p:\tM\rightarrow M$ with branch set $\Sigma$. Thus $M$ is the quotient of $\tM$ by an involution $\tau:\tM\rightarrow \tM$. The additional data in the representation $\chi$ is equivalent to a class in $H^{1}(\tM)^{-}$, the $-1$ eigenspace of the action of $\tau$ on $H^{1}(\tM, \bR)$,  and our assumption is that this is nonzero. This is standard algebraic topology. To see the correspondence in one direction, chose any section $\psi$ of $E^{+}$ which vanishes near $\Sigma$ (in the sense of the local vector bundle structure). Then $\theta=d\psi$ can be interpreted as a closed, $\tau$ anti-invariant,  $1$-form on $\tM$ and we take the de Rham cohomology class $[\theta]$. The discussion extends to the case of sections like $\phi$, except that the resulting form $1$-form $\theta$ may not be smooth. From this point of view the existence and uniqueness of the harmonic section $\phi$ is a version of the Hodge Theorem for harmonic $1$-forms  on $\tM$ with the singular Riemannian metric $p^{*}(g)$.

\item The classical case is when $n=2$ and $M$ is a Riemann surface. Then $\Sigma$ is a finite set and $p:\tM\rightarrow M$ is a holomorphic double branched cover of the familiar kind.  The real $1$-form $\theta=d\phi$ is the real part of a holomorphic $1$-form $\Theta$ on $\tM$ and the square $\Theta^{\otimes 2}$ can be viewed as a meromorphic quadratic differential on $M$.  The condition that $\bA=0$ is just the condition that $\Theta^{\otimes 2}$ is a holomorphic quadratic differential and it then necessarily has zeros at the points of $\Sigma$. Conversely, starting with a holomorphic quadratic differential $\Psi$  on $M$ with simple zeros we can form a double cover $\tM$ on which the square root $\sqrt{\Psi}$ is well-defined. Our main result (Theorem 1), in this case, reduces to a well-known fact that the holomorphic quadratic differentials on $M$ are locally parametrised by the cohomology class of the real part of $\sqrt{\Psi}$  in $H^{1}(\tM)^{-}$.

\item 
The author's main motivation for studying this problem comes from a nonlinear version  developed in \cite{kn:D2} involving \lq\lq branched maximal sections'' which we sketch here and which will be treated at length in a subsequent paper.   An $n$-dimensional submanifold of the indefinite space $\bR^{n,m}$ is called a {\it maximal positive submanifold} if its tangent space at each point is a maximal positive subspace for the indefinite form and if it satisfies the Euler-Lagrange equation associated to the volume functional.
Let $\Gamma_{n,m}$ be the affine extension of the indefinite orthogonal group $O(n,m)$
$$ 0\rightarrow \bR^{n+m}\rightarrow \Gamma_{n,m}\rightarrow O(n,m)\rightarrow 1 . $$ A vector $v\in \bR^{n,m}$ with $v^{2}=-2$ defines a transformation in $O(n,m)$:
$$  w\mapsto w+ (v.w)v,  $$
and we say that an element of $\Gamma_{n,m}$ is a {\it reflection} if it maps to a transformation of this kind in $O(n,m)$. 
Suppose now that $\chi':\pi_{1}(M\setminus \Sigma)\rightarrow \Gamma_{n,m}$ is a homomorphism which takes each linking loop $\delta$ to  a reflection and  form a flat bundle $E_{+}\rightarrow M\setminus L$ with structure group $\Gamma_{n,m}$ and fibre $\bR^{n+m}$. We consider sections $\phi$ of $E_{+}$ which are locally given by maps into $\bR^{n,m}$ with image a maximal positive submanifold. Around a point of $\Sigma$ we require the image to be a  a branched submanifold, modelled transverse to $\Sigma$ on the graph of a 2-valued function ${\rm Re} (z^{\sthreetwo})$. We call such sections $\phi$ {\it branched maximal sections}. The problem then is to develop a deformation theory for such branched maximal sections with respect to small changes in the representation $\chi'$. The linearisation of the maximal condition is an operator of Laplace type  and the condition that the model around $\Sigma$ is preserved is a variant of the condition $\bA=0$ studied in this paper.

\item Another motivation from a nonlinear problem comes from the study of special Lagrangian submanifolds. Let $(X,\omega,\Omega)$ be a Calabi-Yau manifold of complex dimension $n$; so $\omega$ is a K\"ahler form and $\Omega$ is a non-vanishing holomorphic $n$-form. Suppose that $M\subset X$ is a special Lagragian submanifold, {\it i.e.} the restrictions of $\omega$ and the real part of $\Omega$ to $M$ vanish. For this discussion there is no loss in supposing that $X$ is a neighbourhood of the zero section in the cotangent bundle $T^{*}M$ and $\omega$ is the standard symplectic form. Then deformations of $M\subset X$ as a Lagrangian submanifold correspond to the graphs of closed $1$-forms on $M$ and so, locally on $M$,  to the derivatives of functions. The linearisation of the special Lagrangian condition is the Laplace equation on functions on $M$, with the induced metric $g$.  If $\Sigma$ and $\chi$ are as above and if $\bA(\Sigma, g, \chi) =0$ it seems likely that there is a $1$-parameter family of (immersed) special Lagrangians $\iota_{t}:\tM\rightarrow X$ which collapse as $t\rightarrow 0$ to  the double cover map with image $M$.

\item Questions of a similar character to the problem we consider here arise in recent work on non-compactness phenomena for various equations in gauge theory over 3-manifolds and 4-manifolds.  This development began with the work of Taubes \cite{kn:T} and has been studied further by a number of other authors,  for example \cite{kn:HW},  \cite{kn:Tak}.  It seems possible that the methods in this paper could have useful applications to these other equations.
\end{enumerate}

\

Our proof of Theorem 1 is an application of Nash-Moser theory. We invoke a general result of implicit function theorem type from Hamilton's exposition  \cite{kn:H},  which develops a variant due to Zehnder \cite{kn:Z} of that theory. With the appropriate analytical foundations and set-up this general implicit function theorem gives what we need rather quickly (essentially all that remains to do is the calculation in the first part of Section 5 below). In turn,  the analytical foundations for the study of these branched harmonic functions and their asymptotics could be regarded as known material, for example as part of much more general theory developed by Mazzeo \cite{kn:M1}, Mazzeo and Vertman \cite{kn:M2}, and other authors. Thus this paper is to some  extent expository in nature, bringing together these different ideas,  and there is  choice of the amount of background material to include. The author's choice is to treat the Nash-Moser theory as a black box, merely quoting the result needed,  but we attempt to give a largely self-contained account of the analytical foundations and set-up for our problem;  partly with an eye towards further developments of the kind sketched above.

\

This work was partially supported by the Simons Foundation through the Simons Collaboration {\it Special holonomy in geometry, analysis and physics}.
The author is grateful to Rafe Mazzeo and Curt McMullen for discussions about this work. 

\section{Analysis foundations}
\subsection{The flat model}
In this subsection we review some results of Schauder type for the flat model of the pair $(M,\Sigma)$. We work on the product $\bR^{n}=\bC \times \bR^{n-2}$ and we will often represent a point in $ \bR^{n}$ by  $(z,t)$ and write $z=re^{i\theta}$.  Let $V$ be the flat vector bundle   over $\bC^{*}\times \bR^{n-2}$  with fibre $\bR$ and holonomy $-1$.   It will be convenient to regard sections of $V$ as being defined over all of $\bR^{n}$,  taking value zero on the singular set
$\{0\}\times \bR^{n-2}$. In the usual way, we can think of sections of $V$ as multivalued functions on $\bR^{n}$: for example the expression ${\rm Re} z^{\shalf}$ defines a section of $V$. From another point of view we can take the branched cover $(w,t)\mapsto (w^{2}, t)$ and identify sections of $V$ with functions  which are odd with respect to the involution
$(w,t)\mapsto (-w,t)$.  Let $\Delta$ be the Laplace operator, acting on sections of $V$ (where we use the \lq\lq analysts'' sign convention). The main focus in this section is on the inverse operator $\uG$. Let $\cL^{2}_{1}$ be the Hilbert space obtained as the completion of the compactly supported sections under the norm $\Vert \nabla u\Vert_{L^{2}}$.  If $\rho$ is an $L^{2}$-section of $V$ with compact support the standard argument using the Riesz representation theorem and Sobolev inequality shows that there is a unique solution $u$ in $\cL^{2}_{1}$ to the equation $\Delta u =\rho$, understood in the weak sense, and we define  $\uG (\rho)  =u$. One way of seeing this is to change the branched covering map to $(w,t)\mapsto (w^{2}/\vert w\vert, t)$. Then the pull-back of the Euclidean metric is uniformly equivalent to a Euclidean metric and we can fit into the standard theory of uniformly elliptic equations with bounded leading co-efficients, restricted to functions which are odd under change of sign of $w$.

 It is also standard that $\uG$ is defined by a Green's function $G(p,p')$:
  $$  \uG(\rho)(p)=\int_{\bR^{n}} G(p,p') dp' . $$
   To be completely precise,  $G$ is a section of a suitable flat bundle but we suppress that in our notation. There is an explicit representation of $G$ in terms of Bessel functions but the properties we need can be summarised as follows.

\begin{enumerate}
\item $G(p,p')$ is symmetric in $p,p'$ and invariant under translations in the $\bR^{n-2}$ factor and the rotations $S^{1}\times O(n-2)$.    
\item $G$ has scaling homogeneity $(2-n)$:
$$  G(\lambda p,\lambda p') = \lambda^{2-n} G(p,p')$$ for $\lambda>0$. 
\item $G$ is smooth on the complement of the diagonal in $\left( \bC^{*}\times \bR^{n-2}\right)^{2}$ and has a standard pole on the diagonal (i.e. equal to the Newton kernel plus a smooth function).
\item For a point $p=(z,t)$ in the open unit ball $B_{1}$ and a point $p'=(z',t')$ outside the ball there is convergent series representation
\begin{equation}   G(p,p')= \real \left( \sum_{k,\nu\geq 0} a_{k,\nu}(t,p') e^{i (\nu+\shalf) \theta} r^{\nu+\shalf + 2k}\right). \end{equation}
This has all the good properties one could hope for. In particular, as $p$ ranges over a compact subset of $B_{1}$  and $p'$ in a compact subset of the complement of that ball,  the derivatives satisfy $\vert \frac{\partial G}{\partial r} \vert \leq C r^{-\shalf}$,  $\vert \frac{\partial G}{\partial \theta}\vert \leq C r^{\shalf}, \vert \frac{\partial G}{\partial t_{i}}\vert \leq C$, and similarly for higher derivatives. 
\end{enumerate}

\

(Throughout this paper we use the standard convention, writing $C$ for a constant which varies from line to line.)

\

{\bf Remark} 

We can decompose $G$ into Fourier components:
$$  G(p,p')=\sum_{\nu\geq 0} G_{\nu} \cos ((\nu+\shalf)\theta), $$
where each $G_{\nu}$ is a function of $r,r',t,t'$. The first component $G_{0}$ has a simple explicit description. Consider a section $f(r,t)$ of the form
$f=r^{-\shalf} g(r,t) e^{i\theta/2}$. Then one finds that 
$$  \Delta f= r^{-\shalf} (\Delta_{r,t} g) e^{i\theta/2}, $$
where $\Delta_{r,t}$ denotes the usual Laplacian in the half space $\bR^{n-1}_{+}$ with co-ordinates $(r,t)$ and $r>0$. Let $K((r,t),(r',t'))$ be the Green's function for this half-space with Dirichlet boundary conditions. It follows from the identity above that
\begin{equation} G_{0}((r,t), (r',t))= (rr')^{\shalf} K((r,t), (r',t')) . \end{equation}
There is a well-known explicit formula for  $K$ obtained by reflection from the Green's function on $\bR^{n-1}$. 

\

\

Fix an exponent $\alpha\in (0,\shalf)$. We define the H\"older norm on sections of $V$ to be
$$  \Vert s \Vert_{,\alpha} = \sup \frac{ \vert s(p)-s(p')\vert}{\vert p-p'\vert^{\alpha}}, $$ where the supremum is taken over pairs $p=(z,t), p'=(z',t')$ with
$$  \vert p-p'\vert \leq \shalf \ {\rm min} (\vert z\vert, \vert z'\vert). $$
This restriction means that there is no ambiguity in defining $\vert s(p)-s(p')\vert $---we compare using parallel transport along the line segment from $p$ to $p'$. By considering parallel transport around a polygon one sees that
\begin{equation}  \vert s(z,t)\vert \leq C \Vert s \Vert_{,\alpha} \vert z \vert^{\alpha}. \end{equation}
 For a section $s$ which is only defined over a ball we define the $C^{,\alpha}$ norm in the obvious way.

Let $\ucT$ be the set of $n$ commuting vector fields $r \frac{\partial}{\partial r}, \frac{\partial}{\partial \theta}, \frac{\partial}{\partial t_{i}}$ and for $k\geq 1$ let $\ucT_{k}$ be the set of differential operators given by monomials of degree $k$ in these vector fields. So for example if $n=3$ the set $\ucT_{2}$ consists of the six differential operators
$$  \frac{\partial^{2}}{\partial t^{2}},  \frac{\partial^{ 2}}{\partial \theta^{2}}, \frac{\partial^{2}}{\partial t\partial \theta}, r \frac{\partial^{2}}{\partial t\partial r},  r\frac{\partial^{2}}{\partial \theta\partial r}, r^{2} \frac{\partial^{2}}{\partial r^{2}}. $$
We define the $\cD^{k,\alpha}$ norm of a section to be
$$   \Vert s \Vert_{\cD^{k,\alpha}}= \max_{0\leq j\leq k, D\in \ucT_{j} } \Vert D s \Vert_{, \alpha}, $$
with the same remark as above in case the section is only defined over a  ball. Let $\cT$ be the set of vector fields on $\bR^{n}$ which are tangent to the singular set $\{0\}\times \bR^{n-2}$ and let $\cT_{k}$ be the set of differential operators which are sums of products of at most $k$ elements of $\cT$.   So $\ucT_{k}$ is a subset of $\cT_{k}$ and it is easy to check, using (3),  that for any $D\in \cT_{k}$ there is a constant $C_{D}$ such that for sections defined over a ball 
$$  \Vert D s \Vert_{,\alpha}\leq C_{D} \Vert s \Vert_{\cD^{k,\alpha}}. $$

\

The analogue of the Schauder estimates that we need is given by the following.
\begin{prop} For  $\rho$ of compact support in  $B_{1}$ 
$$   \Vert (\uG\rho)\vert_{B_{1}}\Vert_{\cD^{2,\alpha}} \leq  C  \Vert \rho\Vert_{, \alpha}. $$
\end{prop}

{\bf Remarks}
\begin{enumerate}
\item We take it as implicit in the statement  that the right side is finite: i.e. that $\rho$ is $\alpha$-H\"older continuous, and we will use the same convention in other statements below. 
\item This can be proved using estimates for the integral representation of $\uG \rho$, just as in the usual Schauder theory. We omit these here, in part because there the author has written an  exposition of similar estimates in \cite{kn:D1}. In fact the situation here is in some ways simpler than in \cite{kn:D1} because the sections vanish on the singular set, which means that there are fewer cases to consider. We establish some similar estimates related to asymptotics in Proposition 2  below. 

\item The result can be strengthened in that for the terms appearing in the ${\cD}_{2,\alpha}$ norm involving operators of order $0$ or $1$ one can replace the exponent $\alpha$ by $\shalf$. We will make use of this observation in the  proof of Lemma 1 below.
\end{enumerate}

\

The Laplace operator $\Delta$ commutes with $\frac{\partial}{\partial t_{i}}, \frac{\partial}{\partial \theta}$ and satisfies the commutation formula
$$  [\Delta, r\frac{\partial}{\partial r}+ \sum t_{i} \frac{\partial}{\partial t_{i}}] =2 \Delta. $$
Using these facts we can immediately extend Proposition 1 to higher derivatives: for each $k$ there is a $C_{k}$ such that
\begin{equation} \  \Vert \uG\rho\vert_{B_{1}}\Vert_{\cD^{k+2,\alpha}} \leq C_{k} \Vert \rho\Vert_{\cD^{k,\alpha}}.
\end{equation}

\

We now consider the asymptotic description of $\uG (\rho)$ near the singular set, which is the cental topic of this paper.  In terms of the series representation (1) we take  the leading term   $a_{0,0}(t,p')$. Clearly this is defined for all $p'\in R^{n}\setminus (0,t)$ and is invariant under translations in $\bR^{n-2}$ so we can write $$a_{0,0}(t, p')= H(t'-t, z'), $$ with a $\bC$-valued function  $H$ defined on $\bR^{n}\setminus\{0\}$.  This depends only on the Fourier component $G_{0}$ and we obtain from (2) an explicit formula:

\begin{equation}
H(t,z)= \kappa_{n} \frac{z^{\shalf}}{R^{n-1}}, \end{equation}
where $R= (\vert z\vert^{2}+ \vert t\vert^{2} )^{\shalf}$ and the constants are $\kappa_{3}= \pi^{-1}$ and $\kappa_{n}= 2(n-3)/{\rm Vol}(S^{n-2})$ for $n>3$. 

This function $H$ has homogeneity $\sthreetwo-n$:
$$  H(\lambda \tau, \lambda z')= \lambda^{\sthreetwo-n} H(\tau, z'). $$
Similarly, if $H_{i}(\tau,z'), H_{ij}(\tau,z')$ are the derivatives
$$ H_{i}= \frac{\partial H(\tau, z')}{\partial \tau_{i}}\ ,\  H_{ij}= \frac{\partial^{2} H(\tau, z')}{\partial \tau_{i}\partial \tau_{j}}$$ then 
$  H_{i}$ has homogeneity $\shalf-n$ and $H_{ij}$ has homogeneity $-\shalf-n$.
Clearly $H,H_{i}$ and $ H_{ij} $ are bounded on the unit sphere in $\bR^{n}$. We define an integral operator $\uH$ taking compactly supported sections of $V$ to functions on $\bR^{n-2}$ by
$$  \uH(\rho)(t) =\int_{\bR^{n}} H(t'-t, z') \rho(p') dp' . $$
Then
$$    \frac{\partial}{\partial t_{i}} \left[ \uH(\rho) (t)\right] = \int_{\bR^{n}} H_{i}(t'-t,p') \rho(p')dp', $$
and similarly for the second derivatives. 

We use the usual H\"older spaces $C^{r,\beta}$ of functions on $\bR^{n-2}$, with norm defined by the sum of the $L^{\infty}$ norm and the $\beta$-seminorm of the  derivatives of order $r$.
\begin{prop}\begin{enumerate}
\item If $\rho$ has compact support in $B_{1}$ then 

$$\Vert \uH(\rho)\Vert_{C^{1,\alpha+\shalf }}\leq C \Vert \rho\Vert_{,\alpha}. $$
\item With $A= \uH(\rho)$, the section $\uG\rho$ has asymptotic behaviour
$$  \uG(\rho)( z,t) = \real \left( A(t) z^{\shalf}\right) + E(z,t) $$
with $\vert E(z,t)\vert \leq C r^{\sthreetwo} \Vert \rho\Vert_{,\alpha}$. 
\end{enumerate}
\end{prop}

This is  proved by routine estimates for the integral operators. In the first part we just consider, for simplicity, the leading term in the $C^{1,\alpha+\thalf}$-norm of $\uH(\rho)$, i.e the $C^{,\alpha+\shalf}$ norm of the derivative. The estimate for the lower term $\Vert \uH(\rho)\Vert_{L^{\infty}}$ is easier. The estimate stated for the leading term is scale invariant so the hypothesis that $\rho$ is supported in the unit ball is irrelevant. It suffices to estimate 
$$I= \int_{\bR^{n}} H_{1}(t, p') \rho(p') dp' - \int_{\bR^{n}}H_{1}(0, p')\rho(p') dp'. $$ 
We want to show that $$\vert I\vert \leq C \vert t\vert^{\alpha+\shalf} \Vert \rho\Vert_{,\alpha}$$ for some constant $C$.

Write $\delta= \vert t\vert$. As  in the usual Schauder theory we consider separately the contributions to the integrals from the regions $\vert p'\vert\geq 2 \delta$ and $\vert p'\vert \leq 2 \delta$. For the first, we write the contribution to $I$ as
$$  \int_{\vert p'\vert\geq 2 \delta} \left( H_{1}(t, p')- H_{1}(0, p')\right) \rho(p') dp'. $$
For simplicity of notation suppose that the vector $t\in \bR^{n-2}$ lies  on the jth coordinate axis (it will be obvious to the reader how to remove this assumption).  Then by the mean value theorem
$$\vert H_{1}(t,p')- H_{1}(0, p')\vert = \delta \vert H_{j1}(\tau, p')\vert, $$
for some $\tau$ with $\vert \tau\vert \leq \delta$, where $H_{j1}$ is the second derivative, as above.   Using the homogeneity of $H_{j1}$ we get 
$$ \vert H_{1}(t,p')- H_{1}(0, p')\vert \leq C \delta \vert p'\vert^{-\shalf-n}. $$

On the other hand 
$$  \vert \rho(z',t')\vert \leq C \vert z'\vert^{\alpha} \Vert \rho\Vert_{,\alpha}\leq C \vert p'\vert^{\alpha} \Vert \rho \Vert_{,\alpha}  $$
Integrating we get a bound on this contribution to $I$
$$ C \delta \int_{2 \delta}^{\infty} x^{-\thalf-n} x^{\alpha} x^{n-1}dx= C\delta \int_{2\delta}^{\infty} x^{\alpha- 3/2} dx= C \delta^{\alpha+\thalf}.$$ 
For the contribution to $I$ from the region $\vert p'\vert \leq 2 \delta$ we take the two terms in $I$ separately. We have
$$ \vert H_{1}(0,p')\vert \leq C \vert p'\vert^{1/2-n} , $$
and $\vert \rho(p')\vert \leq C \vert p'\vert^{\alpha} 
$ so the estimate for the first term is now
$$ C \int_{0}^{2\delta} x^{\alpha} x^{\thalf-n} x^{n-1} dx = C \delta^{\alpha+1/2}, $$ and similarly for the second term. This completes our discussion of the first item in Proposition 2.

For the second item we consider
$$  \uG(\rho)(z,0)- {\rm Re}\left( z^{1/2} \uH(\rho)(z,0) \right). $$
Write $\vert z\vert = \delta$ and consider separately the contributions to the integral
$$ \int_{\bR^{n}} G((z,0), p') \rho(p') dp' $$
from the regions $\vert p'\vert\leq  2\delta$ and $\vert p'\vert \geq 2 \delta$. The first contribution is bounded by  
$$C \delta^{\alpha} \Vert \rho\Vert_{,\alpha} \int_{\vert p'\vert\leq 2\delta} \vert G((z,0), p')\vert dp'. $$
By the scaling behaviour of the Green's function this is at most
$$   C \delta^{2+\alpha} \Vert \rho\Vert_{,\alpha} \int_{\vert p'\vert \leq 2} \vert G((1,0), p') \vert dp'. $$
(The finiteness of the integral on the right hand side follows from our statements about the Green's function.)
Similarly one finds that the contribution to the integral defining $\uH(\rho)(0)$ from the region $\vert p'\vert \geq 2\delta$ is bounded by $C\delta^{3/2+\alpha} \Vert \rho\vert_{\alpha}$. Thus it suffices to bound
$$  J= \int_{\vert p'\vert \geq 2\delta} \left( G((z,0), p')\rho(p')- {\rm Re} (z^{\shalf} H(0, p')) \right) \rho(p') dp'. $$

The series expansion (1) of the Green's function implies that there is a constant $C$ such that for all $p''$ with $\vert p''\vert =1$ and all $\tilde{z}$ with $\vert \tilde{z}\vert \leq \shalf$ we have
$$  \vert G((\tilde{z},0), p'') - H(0, p'') \tilde{z}^{\shalf} \vert \leq C \vert \tilde{z}\vert^{\sthreetwo} . $$

Writing $p'= \vert p'\vert p''$ and $z= \vert p'\vert \tilde{z}$ in the integral defining $J$ and using the scaling behaviour of $G$ and $H$ we get
$$  \vert J\vert \leq C \vert z\vert^{\sthreetwo} \int_{\vert p'\vert\geq \delta} \vert p'\vert^{\shalf-n} \vert \rho(p') \vert dp' , $$
and our result follows.

Differentiating with respect to $t$ we immediately extend Proposition 2 to higher order estimates
\begin{equation}  \Vert \uH(\rho)\Vert_{C^{1+k, \thalf+\beta}}\leq C \Vert \rho\Vert_{\cD^{k,\alpha}}. \end{equation}

We can also use the same device as  the proof of (4) to get bounds on the derivatives of the error term in item (2) of Proposition  2. For example
\begin{equation}  \vert \nabla E\vert\leq C r^{\shalf} \Vert \rho\Vert_{\cD^{1,\alpha}}. \end{equation}

\

The next term in the series for $G$ is of order $r^{\sthreetwo}$. We set $a_{1,0}(t,p')= L(t'-t, z')$ and define an integral operator
$$  \uL(\rho)(t)= \int_{\bR^{n}} L(t'-t, z') \rho(p') dp'. $$
Then, with $\rho$ as above,  the same arguments  show that 
\begin{equation} \Vert \uL(\rho)\Vert_{k, \shalf + \alpha}\leq C \Vert \rho\Vert_{\cD^{k,\alpha}}. \end{equation}
(There is an  explicit formula like (5) for $L$, but the estimates only depend on general properties such as the homogeneity so the explicit formulae are not particularly relevant.)
Writing $\uL\rho=B$ we have an asymptotic development
\begin{equation}  (G\rho)(z,t)\sim \real \left( A(t) z^{\shalf} + B(t) z^{\sthreetwo}\right)\end{equation}
as $z\rightarrow 0$, with similar bounds on the error term.

\

\subsection{Variable co-efficients} 
 In this subsection we extend the preceding discussion to variable co-efficient operators, mimicking the standard treatment. Recall that $\cT_{2}$ is the set of differential operators which are sums of products of at most two tangential vector fields. 

We say that a differential operator $\tDelta$ over $\bR^{n}$ is {\it admissible} if:
\begin{enumerate} \item  $\tDelta=\Delta+ \cL$ with $\cL$ in ${\cal T}_{2}$. \item $\tDelta$ is elliptic and has divergence form, specifically we assume that
$$\tDelta = W^{-1} \Delta_{g} W$$ where $W$ is a smooth  positive function and $\Delta_{g}$ is the Laplace operator of a smooth Riemannian metric $g$ on $\bR^{n}$.
\item $\tDelta=\Delta$ outside a compact set in $\bR^{n}$.  
\end{enumerate}

These conditions mean that we can apply the Hilbert space theory, as before, so for any $\rho\in L^{2}$ of compact support there is a unique solution
$u$ of the equation $\tDelta u=\rho$ in the Hilbert space $\cL^{2}_{1}$. Throughout this subsection we suppose that $\tDelta$ is an admissible operator.

\begin{prop} Suppose that $\tDelta=\Delta+\cL$ is  equal to $\Delta$ outside a compact subset of $B_{1}$ and that  the coefficients
of $\cL$ are sufficiently small in $C^{2}$. Then if $\rho\in C^{,\alpha}$ has compact support in $B_{1}$ the solution $u$ of $\tDelta u=\rho$ has $\Delta u=\sigma$ with $\sigma$ of compact support in $B_{1}$ and $$\Vert \sigma\Vert_{,\alpha}\leq C \Vert \rho \Vert_{,\alpha}.$$
\end{prop}

To see this we consider the equation $\tDelta(\uG\sigma)=\rho$ for $\sigma$. This is \begin{equation} \sigma+ \cL \uG \sigma = \rho. \end{equation}
By assumption $\cL \uG \sigma$ vanishes outside $B_{1}$ so the equation implies that $\sigma$ is supported in $B_{1}$. Our Schauder estimates of Proposition 1 imply that for any $\sigma$, $$  \Vert \cL \uG \sigma\Vert_{,\alpha} \leq C'_{\cL} \Vert \sigma \Vert_{C^{,\alpha}}, $$
for some constant $C'_{\cL}$ depending on $\cL$. If the coefficients of $\cL$ are sufficiently small we can arrange that $C'_{\cL}<1$.  Then the operator $\cL \uG$ is a contraction, equation (10) has  a solution $\sigma\in C^{,\alpha}$ and, by uniqueness, $\uG\sigma$ agrees with the weak solution.

To formulate a general elliptic regularity result, we
let $\cD^{k,\alpha}_{\loc}$ be the sheaf of sections defined by the $ \cD^{k,\alpha}$ norms, in the standard fashion. Then  we have
\begin{prop} Suppose that $u$ is defined on some open set in $\bR^{n}$ and lies in $\cD^{1,\alpha}_{\loc}$ and $\tDelta u\in \cD^{k,\alpha}_{\loc}$,  then $u\in \cD^{2+k,\alpha}_{\loc}$.\end{prop}

Consider the case $k=0$.
Away from the singular set this is standard elliptic regularity so it suffices to work over a small ball $B_{r}$. Let $\chi$ be a cut off function supported in $B_{r}$ which is the product of cut-of functions of $\vert z\vert $ and $\vert t\vert $ in the  obvious way. This means that  the derivatives of $\chi$ in the $\bC$ factor vanish near the singular set. Then \begin{equation}\tDelta(\chi u) = \chi \tDelta u + 2 \nabla \chi. \nabla u + F u \end{equation}
where $F$ is smooth. The point now is that $\nabla\chi.\nabla u$ only involves the $t$ derivatives of $u$ near the singular set so $\tDelta (\chi u)$ lies in $C^{,\alpha}$. If $r$ is sufficiently small we can suppose that after rescaling to the unit ball  the rescaled operator $\tDelta_{\scaled}$ satisfies the small perturbation assumption of Proposition 3. We also use cut-off functions to  modify $\tDelta_{\scaled}$ to be equal to $\Delta$ outside $B_{1}$. Then Propositions 1 and 3 imply that $\chi u$ is in $\cD^{2,\alpha}$. 

The case of general $k$ follows just as in the standard theory,  using a nested sequence of balls.

\

We will need a stronger result. Let $\cL^{2}_{1,\loc}$ be the sheaf of $L^{2}$ sections with weak derivatives in $L^{2}$. (It is straightforward to show that $\cD^{1,\alpha}_{\loc}$ is contained in $\cL^{2}_{1,\loc}$.)

\begin{prop}

If $u\in \cL^{2}_{1,\loc}$ is a weak solution of the equation $\tDelta u=0$ on an open set $U\subset \bR^{n}$ then $u$ is in $\cD^{k,\alpha}_{\loc}$ for all $k$. More generally, if $u$ is a weak solution of the equation $\tDelta u=\sigma$ for $\sigma\in \cD^{k,\alpha}_{\loc}$ then $u$ is in $\cD^{k+2,\alpha}_{\loc}$.
\end{prop}

To prove this we use the following Lemma.
\begin{lem}
Suppose that $\tDelta$ satisfies the hypotheses of Proposition 3. There is a constant $C$ such that if $\rho\in C^{,\alpha}$ is supported in the annulus $B_{1}\setminus B_{2/3}$ and $v$ is the solution of $\tDelta v=\rho$ given by Proposition 3 then over the interior ball $B_{1/3}$ we have
$$  \Vert \nabla_{t} v\Vert_{\alpha} \leq C \Vert \rho\Vert_{L^{2}}. $$
\end{lem}

Assuming this Lemma, the proof of Proposition 5 follows from an approximation argument. Let $u$ be a weak solution of $\tDelta u=0$. By scaling, we can suppose that $u$ is defined over the unit ball and that $\tDelta$ satisfies the hypotheses of Proposition 3. Let $\chi$ be a cut-off function with derivative supported in the annulus $B_{1}\setminus B_{\stwothirds}$. Then $\tilde{\rho}= \tDelta(\chi u)$ is supported in this annulus and lies in $L^{2}$. Approximate $\tilde{\rho}$ in $L^{2}$ norm by a sequence $\rho_{i}$ in $C^{,\alpha}$, supported in this annulus. It is a simple consequence of the Lemma that $\nabla_{t} u$ is in $C^{,\alpha}$ over the interior ball $B_{\sthird}$. Then the same argument as in the proof of Proposition 4 shows that $u$ is in $\cD^{k,\alpha}_{\loc}$ for all $k$.  The generalisation  to the case when $\tDelta u$ is in $\cD^{k,\alpha}$ is straightforward.

We now turn to the proof of Lemma 1. We assume that $\Vert \rho\Vert_{L^{2}}=1$. Moser iteration gives an $L^{\infty}$ bound: 
$$ \Vert v \Vert_{L^{\infty}}\leq M. $$
For $p,q$ in the interior ball $B_{\shalf}$, write $\delta(p,q)$  for the distance from $\{p,q\}$ to the boundary $\partial B_{\shalf}$ and 
$$ Q(p,q)= \delta(p,q)^{1+\alpha} \vert p-q\vert^{-\alpha}\vert  \nabla_{t} v (p)-\nabla_{t} v(q)\vert. $$
Let $\mu$ be the supremum of $Q(p,q)$ 
where $p,q$ run over pairs with $\vert p-q\vert \leq \squart \delta(p,q)$. Thus a bound on $\mu$ gives a bound on the $C^{,\alpha}$ norm of $\nabla_{t} v$ in the interior ball $B_{\sthird}$, which is what we seek. Let $p,q$ be a pair such that $\vert p- q\vert \leq \squart \delta(p,q)$ and $Q(p,q)\geq \shalf \mu$. We can suppose that $p$ has distance $\delta=\delta(p,q)$ from the boundary. Let $\tilde{B}$ be the ball of radius $\delta/2$ centred on $p$ and rescale to a unit ball $\tilde{B}'$ of radius $1$. There will be various cases,  depending on the relative position of the singular set, but this will not matter.  Let $v'$ be the function on $\tilde{B}'$ corresponding to the restriction of $v$ to $\tilde{B}$. By construction the derivative $\nabla_{t}v'$ satisfies a Holder estimate
$$  \Vert \nabla_{t} v'\Vert_{,\alpha} \leq  \mu 2^{-(1+\alpha)}, $$
and if $p',q'$ are the points in $\tilde{B}'$ corresponding to $p,q$ we have
\begin{equation}   \vert \nabla_{t}v'(p')- \nabla_{t} v'(q')\vert \geq \mu 2^{-(2+\alpha)} \vert p'-q'\vert^{\alpha}. \end{equation}
Let $\chi'$ be a cut-off function on $\tilde{B}'$ of the kind above,  equal to $1$ on the half-sized ball. Thus $\chi'=1$ at $p'$ and $q'$. If $\tDelta'$ is the rescaled operator over $\tilde{B}'$ we have
$$ \Vert \tDelta' (\chi' v')\Vert_{,\alpha}\leq C \mu. $$
We can plainly suppose that $\tDelta'$ is a small perturbation of the flat model, so we can apply Proposition 3 (or an obvious extension of that for balls not centred at the origin). Now we bring in the fact noted in Remark 3 after Proposition 1 that we get a $\shalf$-H\"older estimate for the first derivative. So
$$  \Vert\nabla_{t} (\chi' v') \Vert_{,\shalf}\leq C \mu. $$
Combined with (12), this implies that $\vert p'-q'\vert$ is not small and that
$$  {\rm max} \left( \vert \nabla_{t} v'(p')\vert, \vert \nabla_{t}v'(q')\vert\right)\geq c \mu, $$
for some constant $c>0$. Suppose that $  \vert \nabla_{t} v'(p')\vert\geq c\mu$ and write $\xi=\nabla_{t} v'(p)$.  The H\"older estimate on $\nabla_{t}v'$ implies that $\nabla_{t}v'$ is within $c\mu/2$ of $\xi$ for points within a fixed small distance from $p'$. Integrating the derivative over a suitable segment we find a nearby point $p''$ such that
$$  \vert v'(p')- v'(p'')\vert \geq c''\mu, $$
for some $c''>0$. Since $\Vert v'\Vert_{L^{\infty}}\leq M$ we get $2M\geq c''\mu$ which gives our bound on $\mu$. The argument in the case when $  \vert \nabla_{t} v'(q')\vert\geq
c\mu$ is the same. This completes the proof of Lemma 1 and hence of Proposition 5.

\

To sum up the results of this subsection we introduce sheaves $\cC^{k+2,\alpha}_{\loc}$ defined by sections $u\in \cD^{k+2,\alpha}_{\loc}$ with $\Delta u\in \cD^{k,\alpha}_{\loc}$. By Proposition 4 it is equivalent to suppose that $\tDelta u\in \cD^{k,\alpha}_{\loc}$ for any admissible operator $\tDelta$. It is also equivalent to suppose that $\Delta_{\bC} u$ is in $\cD^{k,\alpha}_{\loc}$, where $\Delta_{\bC}$ is the Laplacian in the $\bC$ factor. We have sheaf maps
$$  A: \cC^{k+2,\alpha}_{\loc}\rightarrow C^{k+1,\alpha+\shalf}_{\bR^{n-2},\loc}\ \ \ ,\ \ \ B: \cC^{k+2,\alpha}_{\loc}\rightarrow C^{k,\alpha+\shalf}_{\bR^{n-2},\loc}, $$
which give the asymptotic behaviour around the singular set.
The sheaf $\cC^{\infty,\alpha}_{\loc}$ will play the role,  in our theory,  of the sheaf of smooth functions in the standard case.

\section{Global theory}

Now let $M$ be a compact $n$-manifold, $\Sigma\subset M$ a codimension 2, co-oriented submanifold and $E$ a flat  real line  bundle as considered in Section 1. We have sheaves $\cD_{\loc}^{k,\alpha}$ of sections of $E$, defined just as in Section 2, using tangential vector fields. The global sections $\cD^{k,\alpha}$ of $E$ are Banach spaces. An explicit norm  on  $\cD^{k,\alpha}$ depends on various choices, such as a metric on $M$, but all such norms are equivalent. We write $\cL^{2}_{1}$ for the space of $L^{2}$ sections with weak derivatives in $L^{2}$. As noted before, $\cD^{1,\alpha}$ is contained in  $\cL^{2}_{1}$.

We define a {\it normal structure} on $\Sigma$ to consist of the following.

\begin{enumerate}\item A normal bundle $N\subset TM\vert_{\Sigma}$, complementary to $T\Sigma$.
\item A Euclidean structure on $N$.
\item A 2-jet along $\Sigma$ of diffeomorphisms from the total space of $N$ to $M$, extending the canonical $1$-jet.
\end{enumerate}

 The group  of diffeomorphisms of $M$ fixing $\Sigma$ pointwise acts transitively on the normal structures.

A Riemannian metric $g$ on $M$ defines a normal structure. We take the normal bundle given by the orthogonal complement of $T\Sigma$ with the induced Euclidean structure and the 2-jet represented by the normal exponential map $\exp_{\Sigma}:N\rightarrow M$. Conversely any normal structure arises from some metric. 

Consider a system of co-ordinate charts $\psi_{a}:U_{a}\rightarrow M$ for $U_{a}$ open in $\bR^{n}=\bC\times \bR^{n-2}$ covering a neighbourhood of $\Sigma$ and compatible with $\Sigma$ in the obvious sense. If $\tDelta$ is a differential operator acting on sections of $E$ these charts define operators  $\tDelta_{a}$ over $U_{a}$. We say that $\tDelta$ is admissible in these charts if the $\tDelta_{a}$ can be extended to admissible operators over $\bR^{n}$,  in the sense of the previous section. Suppose that we have a normal structure on $\Sigma$. If  we are take an open cover of $\Sigma$ with corresponding local trivialisations of the normal bundle $N$ and a diffeomorphism $\psi$ from a neighbourhood of the zero section of $N$ to $M$ representing the 2-jet,  we get a  system of co-ordinate charts as above. We say that such a system of charts is adapted to the normal structure.

Now let $g$ be a Riemannian metric on $M$,  $\Delta_{g}$ the Laplace operator and $\mu=\mu_{g,\Sigma}$ the mean curvature of $\Sigma$, a section of the normal bundle.  Let $W$ be a smooth positive function on $M$, equal to $1$ on $\Sigma$ and with gradient vector field $\nabla W$ equal to $-\shalf\mu$ on $\Sigma$. Define an operator 
$$\tDelta= W^{-1} \Delta_{g} W,  $$ 
acting on sections of $E$. 

\begin{prop}
The operator $\tDelta$ is admissible in any system of charts adapted to the normal structure defined by $g$.
\end{prop}

Recall that the Laplacian $\Delta_{g}$ can be expressed in terms of a local orthornomal frame of vector fields $X_{i}$ as
\begin{equation}   \Delta_{g}= \sum X_{i}^{2} + \dive(X_{i}) X_{i}. \end{equation}

Consider a chart adapted to the normal structure, so we have an open set
$U\subset \bC\times\bR^{n-2}$ and maps
$$ \phi:  U \rightarrow N , \psi:N\rightarrow M, $$
where $\psi$ is equal up to second order to the normal exponential map. Write
$z=x_{1}+ ix_{2}$ for the standard co-ordinate on $\bC$ and $Z_{1}, Z_{2}$ for the vector fields on the image of the chart given by pushing forward
$\frac{\partial}{\partial x_{1}}, \frac{\partial}{\partial x_{2}}$ under
$\psi\circ \phi$. So $Z_{1}, Z_{2}$ restricted to $\Sigma$ give an orthonormal frame for the normal bundle. We claim that on $\Sigma$ 
\begin{equation}  \nabla_{Z_{i}} Z_{j}=0, \end{equation}
(using the Levi-Civita covariant derivative $\nabla$ of the metric $g$). In fact this equation (14) gives another characterisation of the normal structure defined by $g$. To verify the claim we  can assume that $\psi$ is the normal exponential map. Let $\gamma(t)$ be the geodesic emanating from a point $p\in \Sigma$ with initial tangent vector $Z_{1}$. Then along $\gamma$ the vector field $Z_{1}$ is the velocity vector of $\gamma$, which is covariant constant by the geodesic equation. Along $\gamma$ the vector field $Z_{2}$ is $t^{-1} V$ where $V$ is the solution of the Jacobi equation with initial condition $V(0)=0, V'(0)= Z_{2}(p)$. The fact that $\nabla_{Z_{1}} Z_{2}=0$ at $p$ follows from the standard discussion of the small $t$ behaviour  of solutions of the Jacobi equation.  

A consequence of (14) is that $\langle Z_{i}, Z_{j}\rangle = \delta_{ij}+O(r^{2})$ where $r$ is the distance to $\Sigma$. So we can choose orthonormal vector fields $\uZ_{1}, \uZ_{2}$ in the span of $Z_{1}, Z_{2}$ with $\uZ_{i}=Z_{i}+O(r^{2})$. Then on $\Sigma $ we have $\nabla_{\uZ_{i}}\uZ_{j}=0$.   Extend this orthonormal pair to an orthonormal frame $\uZ_{1}, \uZ_{2}, Y_{1}, \dots , Y_{n-2}$: thus the vector fields $Y_{i}$ are tangent to $\Sigma$. We have 
$$\Delta_{g}= \Delta_{1} + \Delta_{2} $$
where $$ \Delta_{1}= \uZ_{1}^{2} + \uZ_{2}^{2} + \dive(\uZ_{1}) \uZ_{1} + \dive(\uZ_{2})\uZ_{2}$$ and $$\Delta_{2}= \sum Y_{i}^{2} + \dive(Y_{i}) Y_{i}. $$
Since the vector fields $Y_{i}$ are tangent to $\Sigma$ the operator $\Delta_{2}$ pulls back over $U$ to an operator in ${\cal T}_{2}$. By construction $Z_{1}^{2}+Z_{2}^{2}$ pulls back to the operator $\Delta_{\bC}$.  Since $\uZ_{1}= Z_{1}+O(r^{2})$ we can write
$$  \uZ_{1}= Z_{1} + \sum f_{a} g_{a} \zeta_{a}$$ for vector fields $\zeta_{a}$ and functions $f_{a}, g_{a}$ vanishing on $\Sigma$. Expanding out,  one finds that $\uZ_{1}^{2}- Z_{1}^{2} $ is in ${\cal T}_{2}$ (i.e. a sum of composites of vector fields vanishing on $\Sigma$).  Combined with same argument for $\uZ_{2}$, this shows  that the pull back of $\uZ_{1}^{2}+ \uZ_{2}^{2}$ is equal to $\Delta_{\bC}$,  up to elements of ${\cal T}_{2}$. 

It remains then to examine the term
\begin{equation} (\dive \uZ_{1}) \uZ_{1}+ (\dive \uZ_{2}) \uZ_{2}. \end{equation}
We will formulate the statement in general codimension.

\begin{lem} Let $S\subset M$ be a codimension-$p$ submanifold of a Riemannian manifold and $\uZ_{1},\dots \uZ_{p}$ be vector fields on $M$ which restrict to an orthonormal frame for the normal bundle of $S$  and which satisfy $\nabla_{\uZ_{i}} \uZ_{j}=0$ on $S$. Then  the vector field
$ \shalf \sum \dive(\uZ_{i}) \uZ_{i}$ restricts to the mean curvature of $S\subset M$.\end{lem}

We leave the proof as an exercise for the reader.

Given this Lemma, we see that if $\Sigma$ is not a minimal submanifold the operator $\Delta_{g}$ is {\it not} an admissible operator with respect to these coordinate charts.  This is the point of the function $W$. We have
\begin{equation}  W^{-1} \Delta_{g} W (u)= \Delta_{g}  u + 2 W^{-1}\nabla W.\nabla u + ( W^{-1} \Delta_{g}( W)) u\end{equation}
By the choice of $W$, the first order term cancels (15), up to a vector field vanishing on $\Sigma$, so $\tDelta$ is admissible.

We define sheaves 
\begin{equation} \cC_{\loc}^{k+2,\alpha} = \{ u \in \cD_{\loc}^{k+2,\alpha}: \tDelta u \in \cD_{\loc}^{k,\alpha}\}. \end{equation}
While this definition involves $\tDelta$, which depends on the metric $g$ and the choice of a function $W$, it follows from the discussion in Section 3 that the sheaves $\cC_{\loc}^{k+2,\alpha}$ depend only on the normal structure. As usual, we have Banach spaces $\cC^{k,\alpha}$ of global sections.

\begin{prop}
Given $g,W$ as above the operator $\tDelta: \cC^{k+2,\alpha}\rightarrow \cD^{k,\alpha}$ is an isomorphism of Banach spaces.
\end{prop}
 
It is immediate from the definition that $\tDelta$ defines a bounded operator between these spaces.  To see that it an isomorphism we go through the Hilbert space theory. This implies,  in a standard way,  that for any $\rho\in \cD^{k,\alpha}$ there is a unique weak $\cL^{2}_{1}$ solution to the equation $\tDelta u=\rho$; then our regularity result (Proposition 5) shows that this is in $\cC^{k+2,\alpha}$.

Of course $\tDelta$ is essentially equivalent to $\Delta$,  so we have that
$$\Delta_{g}:  W \cC^{k+2,\alpha}\rightarrow W \cD^{k,\alpha}$$ is an isomorphism. In fact $W\cD^{k,\alpha}$ is equal to $\cD^{k,\alpha}$.  The point is that, unlike $\cD^{k,\alpha}$,  the function space $\cC^{k+2,\alpha}$ is not preserved by multiplication by smooth functions on $M$.

\

We now have a global version of the asymptotic expansion of sections around $\Sigma$. Given a normal structure on $\Sigma$, fix a compatible diffeomorphism from a neighbourhood 
 of the zero section in $N$ to a neighbourhood $U$ of $\Sigma$ in $M$. The inverse is a map
$  \zeta:U\rightarrow U'\subset N$. Using the co-orientation of $\Sigma\subset M$ we can regard $N$ as a complex line bundle over $\Sigma$. 
If $\sigma_{1}$ is a section of the dual bundle $N^{-1}$ we get a complex valued function $\langle \sigma_{1} , \zeta\rangle$ on $U$. We will just denote this function by $\sigma_{1}\zeta$.   More generally, the flat real line bundle $E$ defines a square root $N^{-\shalf}$ of $N^{-1}$ and if $p$ is an integer or half-integer and $\sigma_{p}$ is a section of $N^{-p}$ we can define  $\sigma_{p}\zeta^{p}$.
 This  is a complex-valued function on $U\setminus \Sigma$ if $p$ is an integer and a section of the complexified bundle $E^{\bC}$ over $U\setminus \Sigma$ if $p$ is not an integer.  Then we have real valued functions and sections of $E$ given by  ${\rm Re} (\sigma_{p}\zeta^{p})$.  
Changing the choice of the diffeomorphism, for the fixed normal structure, changes $\sigma_{p}\zeta^{p}$ by $O(r^{p+2})$ so, for example,  a statement that
$$   u \sim \real (\sigma \zeta^{1/2}+ \tau \zeta^{3/2}) $$
is independent of this choice of diffeomorphism. It will be convenient to suppose that the $\sigma_{p}\zeta^{p}$ are extended (in some arbitrary way) outside $U$. 

\

With this notation in place, we have maps
\begin{equation} A: \cC^{k+2,\alpha}\rightarrow C^{k+1,\alpha+\shalf}(N^{-1/2}), \end{equation}
\begin{equation}B:\cC^{k+2,\alpha}\rightarrow C^{k,\alpha+\shalf}(N^{-3/2}), \end{equation}
such that for any $u\in \cC^{k+2,\alpha}$ we have
$$  u\sim \real(A(u) \zeta^{\shalf} + B(u) \zeta^{\sthreetwo}),  $$
and we have estimates in the manner of Proposition 2 (2) and equation (5) for the error term.

\section{Harmonic section of the affine bundle and Nash-Moser theory}

Recall that we have a flat affine bundle $E_{+}$ lifting $E$, defined by a non-zero class $h$ in a cohomology group which we just denote by $H^{1}$. The space of sections $\Gamma(E_{+})$ is an affine space modelled on the vector space $\Gamma(E)$. The Laplace operator $\Delta_{g}$ is defined on sections of $E_{+}$, mapping to sections of $E$. In a small neighbourhood of each connected component of $\Sigma$ there is a canonical isomorphism of $E_{+}$ with $E$. We choose the function $W$ to be equal to $1$  outside these small neighbourhoods and define $W:\Gamma(E^{+})\rightarrow \Gamma(E^{+})$ to be given by multiplication by $W$ inside these neighbourhoods, using the identification with $E$,  and the identity elsewhere. So we have an operator $$\tDelta= W^{-1}\Delta_{g} W  :\Gamma(E^{+})\rightarrow \Gamma(E)$$ and for $\phi_{1}, \phi_{2}\in \Gamma(E_{+})$ the difference $\tDelta(\phi_{1})-\tDelta(\phi_{2})$ is equal to $\tDelta(\phi_{1}-\phi_{2})$ as considered before. We have spaces $\cC^{k,\alpha}(E_{+})$ and $\tDelta$ is an isomorphism from $\cC^{k+2,\alpha}(E_{+})$ to $\cD^{k,\alpha}(E)$. In particular there is a unique section $\shar$ in $\cC^{\infty,\alpha}(E^{+})$ with $\tDelta \shar=0$ and $\phi= W \shar $ is the unique harmonic section of $E_{+}$. Since $E_{+}$ is identified with $E$ near $\Sigma$ the discussion of asymptotics is identical and we have $A=A(\shar)\in C^{\infty}(N^{-\shalf})$ and $B=B(\shar)\in C^{\infty}(N^{-\sthreetwo})$. The $O(r^{\sthreetwo})$ term in the asymptotics of $\phi$ is affected by the mean curvature and we have
\begin{equation}
  \phi= {\rm Re}(A \zeta^{\shalf} + B\zeta^{\sthreetwo}) -\shalf {\rm Re}(A \zeta^{\shalf}) {\rm Re}(\overline{\mu} \zeta) + O(r^{\sfivetwo}), \end{equation}
where $\overline{\mu}\in \Gamma(N^{-1})$ corresponds to $\mu$ under the antilinear isomorphism from $N$ to $N^{-1}$ furnished by the metric.

\

We will now digress to recall some of the main elements of the Nash-Moser theory, following Hamilton's article \cite{kn:H}. A fundamental notion  is that of a {\it tame estimate}. The initial context for this is in the setting of Fr\'echet spaces with an increasing sequence of norms:
$$  \Vert f\Vert_{0}\leq \Vert f\Vert_{1}\leq \Vert f\vert_{2}\dots. $$
A map $S$ from an open set in one such space to another satisfies a tame estimate if there is an $r$ and  constants $C_{m}$ such that
\begin{equation}   \Vert S(f)\Vert_{m}\leq C_{m}\left( 1+\Vert f\Vert_{m+r}\right), \end{equation}
for all sufficiently large $m$ and for all $f$ in the domain of the map (\cite{kn:H} II.1.3.2).  A {\it smooth tame map} is a smooth map from such an open set all of whose derivatives satisfy tame estimates. There is also a notion of a tame Fr\'echet space (\cite{kn:H}, II.2.1.1),  which we do not need to recall here.  The definitions extend to a class of tame Fr\'echet manifolds,  and maps between them.

Returning to our set-up, fix a submanifold $\Sigma_{0}$ with normal structure. Let $\cM$ be the  space of  Riemannian metrics on $M$: this is an open subset in a tame Fr\'echet space.  Let $\cM_{0}\subset \cM$ be the subset of metrics compatible with the normal structure. This imposes algebraic conditions on the $1$-jet of the metric along $\Sigma$ and it is easy to check  that $\cM_{0}$  is a tame Fr\'echet submanifold. It is also easy to define a smooth tame map $w:\cM_{0}\rightarrow C^{\infty}(M)$ such that $W=w(g)$ is a function of the kind considered above, compatible with the mean curvature of $\Sigma_{0}$ in the metric $g$.
 We have a map
$$  \uA: H^{1}\times \cM_{0}\rightarrow C^{\infty}(\Sigma_{0}, N^{-1/2}) . $$
 defined by $A(\shar)$,  where $\shar$ is the section above corresponding to the pair $h,g$ and $w(g)$.  
\begin{prop}
$\uA$ is a smooth tame map. 
\end{prop}

For fixed $g$ the map $\uA$ is linear in $h$ so the $h$ dependence is straightforward and for simplicity we will fix $h$. We can fix a section $\phi_{0}$ of $E_{+}$ which is covariant constant near $\Sigma$ and write $\shar= \phi_{0}+U$ for a section $U$ of $E$.   
Using $W=w(g)$, we have a family of operators $\tDelta_{g}$ acting on a fixed scale of Banach spaces
$$  \tDelta_{g}: \cC^{k+2,\alpha}(E)\rightarrow \cD^{k,\alpha}(E) $$ and the sections $\shar$ correspond to solutions of a family of equations
\begin{equation} \tDelta_{g}(U)= - \rho_{g}, \end{equation} where $\rho_{g}= -\tDelta_{g}\phi_{0}$.

Fix a reference metric $g_{0}\in \cM_{0}$ and write $g=g_{0}+ \gamma$.
Let $\Vert\ \Vert_{k}$ be the usual  $C^{k,\alpha}$ H\"older norm on the space of metric tensors and work in a neighbourhood $\Vert \gamma\Vert_{K}<\epsilon$ for a suitable $K$ and suitable small $\epsilon$.  
We know that there is a unique  solution $U$ to (22)  and we want to show that it satisfies  tame estimates, for some suitable $r$,
\begin{equation} \Vert U\Vert_{\cC^{m,\alpha}}\leq C_{m} \left(1+ \Vert \gamma \Vert_{m+r}\right),\end{equation}  on this neighbourhood in  $ \cM_{0}$.   
Write $\tDelta_{g}=\tDelta_{g_{0}}+ \cL$ and let $\cF$ denote the coefficients of the operator $\cL$---{\it i.e.} a section of the dual of a 2-jet bundle. Using the same notation to denote norms of these sections, we clearly have a tame estimate on $\cF$ in terms of $\gamma$. In fact we have
\begin{equation}  \Vert \cF\Vert_{m}\leq C_{m} \Vert \gamma \Vert_{m+1}.\end{equation}

The proof now is similar to that of Theorem 3.3.1 in \cite{kn:H} but we cannot immediately quote that result because of the special character of our function spaces $\cC^{k,\alpha}, \cD^{k,\alpha}$. Writing $\tDelta_{0}=\tDelta_{g_{0}}$, we can take the $\cC^{k,\alpha}$ norm to be defined by
$$   \Vert u\Vert_{\cC^{k,\alpha}}= \Vert \tDelta_{0} u \Vert_{\cD^{k,\alpha}}. $$ Taking suitable $\epsilon,K$ we can treat $\tDelta_{g}$ as a small perturbation of $\tDelta_{0}$, in the manner of Subsection 2.2 above, and we get
\begin{equation}  \Vert u \Vert_{\cC^{2,\alpha}}\leq C \Vert \tDelta_{g} u \Vert_{,\alpha}\end{equation}
We now follow the usual strategy of differentiating this equation to get estimates for higher derivatives. Let $X_{1}, \dots X_{k}$ be tangential vector fields and $\nabla_{i}=\nabla_{X_{i}}$. Let $T_{i}$ be the commutator $T_{i}=[\tDelta_{g}, \nabla_{i}]$. We have

\begin{equation} [\tDelta_{g}, \nabla_{1}\dots \nabla_{k}] = \nabla_{1}\dots\nabla_{k-1}T_{k}+\dots T_{1}\nabla_{2}\dots \nabla_{k}, \end{equation}
with a sum of $k$ terms on the right hand side.
The operators $T_{i}$ are of second order and it follows from the definitions that they can be expressed in terms of tangential vector fields, so for {\it fixed} $g$ we have

$$   \Vert [\tDelta_{g}, \nabla_{1}\dots\nabla_{k}]u \Vert_{\cD^{k,\alpha}}\leq C(g) \Vert u \Vert_{\cD^{k+2,\alpha}}. $$

A small subtlety now arises which it is easiest to explain by an example in local co-ordinates. Take $z=x_{1}+ i x_{2}$ as a standard co-ordinate transverse to $\Sigma$ as in Section 2 and consider an operator
$$   f(x_{1}, x_{2}, t) \frac{\partial}{\partial x_{1}}, $$
with $f(0,0,t)=0$. We can express this operator as a linear combination of standard tangential vector fields
$$   h_{1} x_{1} \frac{\partial}{\partial x_{1}} + h_{2} x_{2} \frac{\partial}{\partial x_{1}}, $$
{\it i.e.} 
$$  f = x_{1} h_{1}+ x_{2} h_{2}, $$
but we cannot control the $C^{\ ,\alpha}$ norm of $h_{1}$ and $h_{2}$ in terms of the $C^{,\alpha}$ norm of $f$. But we do have such control in terms of the $C^{1,\alpha}$ norm of $f$. Similarly for second order operators but with the loss of two derivatives.  The upshot of this is that, compared to the usual situation, we lose two orders of differentiability in estimating the right hand side of (26) in terms of $\cF$. But this will not matter for the tame estimate.

What we deduce from (26) is an estimate
$$   \Vert \tDelta_{g}(\nabla_{1}\dots\nabla_{k} u) \Vert_{,\alpha} \leq \Vert \nabla_{1}\dots\nabla_{k} \tDelta_{g} u \Vert_{,\alpha} +  C_{k} \left( \Vert \cF\Vert_{3}\Vert u\Vert_{\cD^{k+1,\alpha}} + \dots \Vert \cF\Vert_{k+2} \Vert u\Vert_{\cD^{2,\alpha}}\right). $$

Now apply this to the solution $U$ of the equation $\tDelta_{g}U=\rho_{g}$ and recall that $\rho_{g}$ is supported away from $\Sigma$. It is clear that $$\Vert \rho_{g}\Vert_{\cD^{k,\alpha}}\leq C_{k} \left( \Vert \gamma \Vert_{k}+1\right). $$ 

Write $N_{p}=\Vert u\Vert_{\cC^{p,\alpha}}$ and $\mu_{q}= \Vert \gamma\Vert_{3+q}$. Putting together the estimates above we get:

\begin{equation} N_{k+2}\leq C_{k}\left( 1+ \mu_{k-2}+\mu_{k}N_{2}+ \dots\mu_{1} N_{k+1}\right) \end{equation}
The Holder norms $\Vert \ \Vert_{m}$ satisfy  interpolation inequalities (see \cite{kn:H} Theorem II.2.2.1 and the remark following). In terms of the $\mu_{q}$ these give, for $j<k$:
$$     \mu_{j}\leq c_{j,k} \ \mu_{k}^{j/k} \mu_{0}^{(k-j)/k}. $$
We can suppose the parameters $\epsilon,K$ chosen so that $\mu_{0}\leq 1$ say. Then a simple induction using (27) shows that there are constants $C'_{k}$ such that
$$   N_{k}\leq C'_{k}(1+ \mu_{k}), $$
which is our tame estimate. 
\

Straightforward arguments similar to (\cite{kn:H} II.3.1.1) show that the map taking $g,h$ to the solution $U$ is a smooth tame map to $\cC^{\infty,\alpha}$. Since the linear map $A$ is bounded from $\cC^{k+2,\alpha}$ to  $C^{k+1,\alpha+\thalf}
$ we deduce immediately that $\uA$ is also a smooth tame map.

\

  Let $\cS$ be the tame Fr\'echet manifold of codimension 2 submanifolds of $M$ (\cite{kn:H} III.2.3.7) and ${\rm Diff}(M)$ the tame Fr\'echet Lie group of  diffeomorphisms of $M$ (\cite{kn:H} III.2.3.5).  For our purposes we will only really need to work in a neighbourhood   of $\Sigma_{0}$  in $\cS$ and a neighbourhood  of the identity in ${\rm Diff}(M)$. We have a smooth tame map from ${\rm Diff}(M)$ to $\cS$ which takes a diffeomorphism $f$ to 
  $f(\Sigma_{0})$.  Let $\cM$ be the space of metrics on $M$, as before and define a subset
   $\cZ\subset \cM\times {\rm Diff}(M)$ to consist of pairs $(g,f)$ such that $f$ takes the normal structure of $\Sigma_{0}$ to the normal structure of $f(\Sigma_{0})$ defined by the metric $g$. It is straightforward to check that $\cZ$ is a tame submanifold. We have  a smooth tame map $$\pi: \cZ\rightarrow \cM\times \cS$$
   which takes $(g,f)$ to $(g, f(\Sigma_{0}))$. In fact this is a tame principal bundle with structure group ${\rm Diff}_{0}$: the group of diffeomorphisms of $M$ preserving $\Sigma$ and its normal structure.

\begin{prop} For a suitable neighbourhood $\cU\subset \cM\times \cS$ of $(g_{0}, \Sigma_{0})$ 
we can find a smooth tame map $\Psi: \cU  \rightarrow {\rm Diff}(M)$ such that $\Psi(g,\Sigma)$ maps $\Sigma_{0}$ to $\Sigma$ and the normal structure of $\Sigma_{0}$ to the normal structure of $\Sigma$ defined by the metric $g$.
\end{prop}
 Another way of expressing this is that $(pr_{\cM}, \Psi)$ is a section of $\pi:\cZ\rightarrow \cM\times \cS$ over $\cU$, where
 $pr_{\cM}$ is the projection onto the $\cM$ factor. 
 
 The existence of a map $\Psi$ with these properties follows from general theory, but we can write down explicit constructions.
 We can suppose that $\Sigma$ lies in a fixed tubular neighbourhood of $\Sigma_{0}$ and can be identified with a section of the normal bundle $N_{0}$. Using this it is easy to write down a diffeomorphism $\psi_{0}$ from a neighbourhood of $\Sigma_{0}$ in $M$ to a neighbourhood of $\Sigma$ in $M$. Let $p$ be a point in $\Sigma_{0}$ and $q=\psi_{0}(p)\in \Sigma$. The derivative of $\psi_{0}$ gives a linear map from $N_{0,p}$ to $TM_{q}$ where $N_{0,p}$ is the normal bundle of $\Sigma_{0}$ at $p$. The image of this map will in general be different from $N_{q}\subset TM_{q}$---the normal bundle defined by the metric $g$---but we can assume that the two subspaces are close, so that orthogonal projection in $TM_{q}$ gives a linear isomorphism $\pi: N_{0,p}
 \rightarrow N_{q}$. Now set 
 $$   \tilde{\pi}= \pi \circ (\pi^{*}\pi)^{-1/2}: N_{0,p}\rightarrow N_{q}, $$
 where the adjoint is formed using the metrics $g_{0}$ on $N_{0,p}$ and $g$ on $N_{q}$. Then $\tilde{\pi}$ is an isometry between the normal bundles with the given metrics. In this way we get an isomorphism
 $\psi_{1}:N_{0}\rightarrow N_{\Sigma}$ of Euclidean vector bundles, covering $\psi_{0}:\Sigma_{0}\rightarrow \Sigma$.  
We have normal exponential maps $$\exp_{\Sigma_{0}}:\Gamma(N_{0})\rightarrow M\ \ ,\ \ \exp_{\Sigma}: \Gamma(N_{\Sigma})\rightarrow M, $$
where the first is defined using the metric $g_{0}$ and the second using the metric $g$.
The composite $\exp_{\Sigma}\circ \psi_{1}\circ \exp_{\Sigma_{0}}^{-1}$ defines a diffeomorphism from a neighbourhood of $\Sigma_{0}$ to a neighbourhood of $\Sigma$,  extending $\psi_{0}$ and taking the normal structure of $\Sigma_{0}$ to the normal structure of $\Sigma$ defined by $g$. It is straightforward, using suitable cut-off functions, to extend this diffeomorphism over the whole of the $M$, equal to the identity outside a larger neighbourhood of $\Sigma_{0}$, and this defines a diffeomorphism $\Psi=\Psi(g,\Sigma):M\rightarrow M$ with the required properties. We leave the reader to check that this is a smooth tame map from $\cM\times \cS$ to ${\rm Diff}(M)$. 

We now have a smooth tame map from $\cU\subset \cM\times \cS$ to $\cM_{0}$, taking $(g,\Sigma)$ to $\Psi^{*}(g)$ where 
$\Psi=\Psi(g,\Sigma)$. In other words, we have used the action of the diffeomorphisms to reduce the study of a pair $(g,\Sigma)$ to the case of the fixed submanifold $\Sigma_{0}$ and the metrics in $\cM_{0}$,  compatible with a fixed normal structure. 

\

The upshot  is that we have a smooth tame map 
$$   \bA: \cU \times H^{1} \rightarrow \Gamma(\Sigma_{0}, N_{0}^{-1/2}), $$

taking $(g,\Sigma, h)$ to $\uA(\Psi^{*}(g,\Sigma)(g), h)$. We want to show that, under suitable conditions and with respect to small variations,  the equation $\bA(g,h,\Sigma)=0$ defines $\Sigma$ implicitly as a function of $g$ and $h$.  For this we quote some deep general results from \cite{kn:H}. 
 
 \
 
 We  introduce  terminology which will simplify the statements of the result we need. 
 Suppose that $ W_{0},  W_{1}, W_{2}$ are tame Fr\'echet spaces and $U$ is an open subset in another tame Fr\'echet space. Suppose that we are given  smooth tame maps $$a:U\rightarrow W_{0} \ \ ,\ \lambda:U\times W_{1}\rightarrow W_{2}\ \mu:U\times W_{2}\rightarrow W_{1}, $$ with $\lambda,\mu$ linear in their second arguments.  Thus we can view $\lambda$ and $\mu$ as  families of linear maps parametrised by $U$. Write
 $$ \mu* \lambda:U\times W_{1}\rightarrow W_{1}\ \ \ \ \ , \  \ \ \lambda*\mu:U\times W_{2}\rightarrow W_{2}, $$
 for the maps defined by composition of these families, in the obvious sense. 
 Then we  will say that  $\mu$ is an {\it inverse to $\lambda$ with $a$-quadratic error} if there are smooth tame maps 
$$Q_{1}: U\times W_{0}\times W_{1} \rightarrow W_{1}\ \ \ ,\ \ \ \ Q_{2}: U\times W_{0}\times W_{2}\rightarrow W_{2}, $$ 
each bilinear in the last two arguments and such that
$$   (\mu*\lambda) (\tau,w_{1}) = w_{1} + Q_{1}(\tau, a(\tau), w_{1})\ \ \ \ (\lambda*\mu)(\tau, w_{2})= w_{2} + Q_{2}(\tau, a(\tau), w_{2}) . $$

With this terminology in place,  let $X,Y, V$ be tame F\'rechet spaces, $U\subset X\times Y $ an open subset
and $A:U\rightarrow V$ a smooth tame map. For $\tau\in U$ let $(D^{Y}A)_{\tau}=\lambda_{\tau}$ be the partial derivative of $A$ in the $Y$ factor, at the point $\tau$.  Thus we are in the setting above, with $W_{0}=W_{2}=V$, with $ W_{1}=Y$,  with $a=A$  and with  $U\subset X\times Y$. 

\begin{thm} Suppose that  $(x_{0}, y_{0})\in U$ and $A(x_{0}, y_{0})=0$. Suppose that there is an inverse to $D^{Y}A$ with $A$-quadratic error. Then there are neighbourhoods $B_{X}$ of $x_{0}\in X$ and $B_{Y}$ of $y_{0}\in Y$ so that for each $x\in B_{X}$ there is a unique $y\in B_{Y}$ with  $A(x,y)=0$ and the map taking $x$ to $y$ is a smooth tame map. 
\end{thm}

This is a combination of Theorems 3.3.1 and 3.3.3 in \cite{kn:H} (which give existence and uniqueness respectively).  Note that if we worked instead in a Banach space setting then the hypotheses would imply that, in a suitably small neighbourhood, the partial derivative $D^{Y}A$ is invertible and the result would follow immediately from the standard implicit function theorem. The use of such \lq\lq approximate inverses''  in the Nash-Moser theory seems to go back to Zehnder \cite{kn:Z}. 

\subsection{The connection}

The representation of our problem as a map $\uA$  from the product  $\cU\subset \cM\times \cS$ to a fixed vector space depends on the choice of the  map $\Psi$ which is far from unique,   so   the partial derivative does not have any intrinsic meaning. The invariant geometric picture involves  a vector bundle over 
$\cM\times \cS$ with fibre $\Gamma(\Sigma, N_{\Sigma}^{-\shalf})$ over $(g,\Sigma)$ (where $N_{\Sigma}$ is the normal bundle defined by the metric $g$). Then we have a section $\bA$ of this bundle and a more intrinsic notion is the covariant derivative of $\bA$ with respect to a connection on the bundle. Indeed Hamilton states a version of Theorem 2  in this  setting (\cite{kn:H}, Theorem 3.3.4).  But the proof that the bundle and section has the required properties will involve the construction of  a local trivialisation and this will amount to  much the same as the approach we took above. However the invariant point of view  makes the discussion of the derivative much more transparent, so we will adopt this language now. We leave the reader to check, that all the constructions below take place in the smooth tame category.

Let $\cG_{0}$ be the automorphism group of the oriented Euclidean vector bundle $N_{0}\rightarrow \Sigma_{0}$. Thus there is an exact sequence
$$    1\rightarrow C^{\infty}(\Sigma_{0}, S^{1})\rightarrow \cG_{0}\rightarrow {\rm Diff}(\Sigma_{0})\rightarrow 1. $$
The corresponding Lie algebra sequence is
$$  0\rightarrow C^{\infty}(\Sigma)\rightarrow \Lie(\cG_{0})\rightarrow {\rm Vect}(\Sigma)\rightarrow 0 . $$
 Let $\cP$ be the space of triples $(g,\Sigma, \psi)$ where $\psi_{1}$ is an isomorphism from $N_{0}$ to $N_{\Sigma}$ covering a diffeomorphism $\psi_{0}:\Sigma_{0}\rightarrow \Sigma$, and $N_{\Sigma}$ is the normal bundle defined by  $g$. The group $\cG_{0}$ acts on $\cP$, making $\cP$ a principal $\cG_{0}$ bundle over $\cM\times \cS$. (This is compatible with the previous discussion of the ${\rm Diff}_{0}$ bundle $\cZ\rightarrow \cM\times \cS$ since the action on the normal bundle of $\Sigma_{0}$ gives a homomorphism from ${\rm Diff_{0}}$ to $\cG_{0}$ and $\cP$ can be induced from $\cZ$ via this homomorphism.)  
The group $\cG_{0}$ acts on the sections of $\Gamma(N_{0})$. For our purposes we can assume that it also acts on $\Gamma(N_{0}^{-1/2})$ (since we are ultimately concerned only with small deformations). Then we get our  vector bundle over $\cM\times \cS$ as  the bundle associated  to this action. 
 
 We want to define a connection on the $\cG_{0}$-bundle $\cP\rightarrow \cM\times\cS$. More precisely, we only need a \lq\lq partial connection'' in the $\cS$-factor.  Fix a point $(g,\Sigma,\psi)$ in $\cP$ and let $T_{\cS}\cP$ be the subspace of $T\cP$ projecting to $T\cS\times \{0\}\subset T(\cS\times \cM)$. So  we have an exact sequence
\begin{equation}  0\rightarrow {\rm Lie}(\cG_{0})\rightarrow T_{\cS}\cP\rightarrow T\cS \end{equation}

The space $T_{\cS}\cP$ can be identified with a subspace of the space $J^{1}$ of 1-jets of sections of $TM$ over $\Sigma\subset M$. That is, if $\tilde{v}$ is the extension of such a $1$-jet to a vector field on $M$ we deform $\psi_{1}$ infinitesimally by composing with the diffeomorphisms generated by $\tilde{v}$.  The condition that this deformation lies in $\cP$, to first order, can be expressed in terms of the Killing operator 
 $$K:\Gamma(M, TM)\rightarrow \Gamma(M, s^{2}TM). $$
 Recall that this is the composite of  the covariant derivative $v\mapsto \nabla v$ with the symmetrisation $TM\otimes TM\rightarrow s^{2}(TM)$. The Killing operator induces a linear map from $1$-jets of vector fields along $\Sigma$ to $0$-jets. {\it i.e.} sections of $s^{2}(TM\vert_{\Sigma})$. We have 
$$  s^{2}(TM\vert_{\Sigma})= s^{2}(N)\oplus (N\otimes T\Sigma) \oplus s^{2}(T\Sigma). $$
Let $$K_{N}:J^{1}\rightarrow \Gamma(s^{2}N\oplus (N\otimes T\Sigma))$$
 be the composite of the Killing operator and the projection on to the two factors indicated. Then, differentiating the normal structure constraint, one finds that  $T_{\cS}\cP$ is the kernel of $K_{N}$. Now  the map from $1$-jets to $0$-jets gives
$$  ev:J^{1}\rightarrow \Gamma(N)\oplus \Gamma(T\Sigma), $$
and the Levi-Civita connection of $g$ gives 
$$  \nabla_{N}: J^{1}\rightarrow \Gamma(N^{*}\otimes TM\vert_{\Sigma}). $$
Let $\cH\subset J^{1}$ be the set of jets $v$ such that $ev(v)$ lies in $\Gamma(N)$ and $\nabla_{N}v=0$. Then the reader will easily check that $\cH$ lies in $T_{\cS}\cP$ and is complementary  to the \lq\lq vertical'' subspace,  thus giving a splitting of the sequence (28). This defines our partial connection.

 Using this partial connection, we have an intrinsic partial covariant derivative of the section $\bA$ over $\cM\times \cS$ which we write as
$$ \frac{\delta \bA}{\delta \Sigma}: \Gamma(\Sigma, N_{\Sigma})\rightarrow \Gamma (\Sigma, N^{-1/2}_{\Sigma}). $$
To relate this to our previous discussion, we suppose for simplicity of notation that we can take $\cU=\cM\times \cS$ so we have a map  $\Psi:\cM\times \cS\rightarrow {\rm Diff}(M)$  as considered in the previous subsection. This induces a section of $\cP\rightarrow \cM\times \cS$ and the projection to the vertical subspace defined using our partial connection gives at each point $(g,\Sigma)$ of $\cM\times \cS$ a linear map
$$  \xi(g,\Sigma): T\cS_{\Sigma}\rightarrow {\rm Lie}(\cG_{0}). $$
This is just the familiar Lie algebra valued $1$-form representing a connection, in our infinite-dimensional setting. The partial derivative computed in a local trivialisation, as considered in the previous section, and the covariant derivative defined above are related by the usual formula
\begin{equation}  \frac{\delta \bA}{\delta \Sigma}(v)= \left(D^{\cS}\bA\right)(v) + \xi(v) (\bA) \end{equation}
where the last term involves the Lie algebra action of ${\rm Lie}(\cG_{0})$ on $\Gamma(N_{0}^{-\shalf})$. (Implicit in our notation here is the identification of the data over $\Sigma$ and $\Sigma_{0}$ furnished by $\Psi$.) The point now is that the term $\xi(v)\bA$ is bilinear in $v$ and $\bA$. This means that (just as in \cite{kn:H} Theorem III.3.3.4) an inverse with $\bA$-quadratic error for $\frac{\delta \bA}{\delta\Sigma}$ is equivalent to one for $D^{\cS}\bA$.  (If we have an explicit choice of $\Psi$, such as that described above, it is possible to write down formulae for $\xi(g,\Sigma)$, but these are very complicated and in the end irrelevant.)

\

Recall that  the second term in the asymptotic expansion gives a section
$\bB(\Sigma, g, h)\in \Gamma(\Sigma; N_{\Sigma}^{-3/2})$. If this is nowhere-vanishing on $\Sigma$ the inverse $\bB^{-1}$ is a section of $N^{3/2}$ and multiplication by $\bB^{-1}$ is a bundle map
$$   \bB^{-1}:\Gamma(N^{-1/2}_{\Sigma} )\rightarrow \Gamma(N_{\Sigma}).  $$ We can now state the principal technical result of this paper.
\begin{thm}
If $\bB(\Sigma,g,h)$ does not vanish anywhere on $\Sigma$ then $\sthreetwo \bB^{-1}$  is an inverse with $A$-quadratic error to the covariant derivative $\frac{\delta A}{\delta \Sigma}$.
 \end{thm}  In this statement it is understood that everything os represented in local trivialisations, putting us in the setting of maps on open subsets of fixed Fr\'echet spaces. The proof of this Theorem is given in 5.1 and 5.2 below. 
 Combined with the Hamilton-Nash-Moser-Zehnder result (Theorem 2) this establishes our main result (Theorem 1) stated in the Introduction.

\section{Calculation of the derivative}

Let $g$ be a Riemannian metric on $M$ and $N$ the normal bundle of $\Sigma$.
If $\sigma$ is a section of $N^{1/2}$ over $\Sigma$ we have defined a section
$\sigma \zeta^{-1/2}$ of $E^{\bC}$  over $M\setminus \Sigma$. Up to order
$r^{3/2}$ it is independent of the choice of $\zeta$. A choice of a function $W$ gives an operator $\tDelta$. 
\begin{thm} Given a choice of $W$, there is a smooth function $f_{W}$ on $M$ with $f_{W}$ and $\nabla f_{W}$ vanishing on  $\Sigma$ and the following property. 
For any $\sigma\in \Gamma(\Sigma,N^{1/2})$ there is a unique section $\tQ=\tQ(\sigma)$ of $E_{\bC}$ over $M\setminus \Sigma$ with $\tDelta \tQ=0$ and such that
$$   \tQ= \sigma \zeta^{-1/2} + f_{W}  \tau\zeta^{-1/2} + u_{\sigma}$$ for some $\tau\in \Gamma(\Sigma, N^{\shalf})$  and  $u_{\sigma}$  in  $\cC^{\infty,\alpha}$. We have  
$$  \Vert u_{\sigma}\Vert_{ \cC^{k,\alpha}}+ \Vert \tau\Vert_{C^{k+2,\alpha}}\leq C_{k} \Vert \sigma\Vert_{C^{k+4,\alpha}}. $$

\end{thm}
If this is true for one choice of $W$ it is true for all, since multiplication by a function equal to 1 to first order on $\Sigma$ preserves the spaces $\cC^{\infty,\alpha}$. For the proof  
it is convenient to work with a particular choice of $W$. There is a standard volume form on the total space of the normal bundle $N$, so we get a function $V$ on a neighbourhood of $\Sigma$ in $M$ such that the pull back of the Riemannian volume of $M$under the normal exponential map is $V$ times this volume form on $N$. Then  we set $W=V^{-1/2}$. It is straightforward to check, and will be clear from the calculations below, that  the derivative of this function is given by the mean curvature as required. Let $\kappa$ be the function on $\Sigma$ given by
$$ \kappa= \squart K_{N} - \shalf {\rm Tr}_{N} {\rm Ric} +\sthreequart \vert \mu\vert^{2} $$
where $K_{N}$ is the sectional curvature in the normal plane and ${\rm Tr}_{N} {\rm Ric}$ is the trace of the Ricci curvature of $M$ restricted to $N$. So in terms of an orthonormal basis $e_{1}, e_{2}$ for $N$, 
$$  {\rm Tr}_{N} {\rm Ric}=  \langle {\rm Ric}\  e_{1}, e_{1}\rangle + \langle {\rm Ric}\  e_{2}, e_{2}\rangle.$$ 
 
 \begin{prop} Write $f=\sigma\zeta^{-\shalf}$.
  \begin{enumerate}\item For any smooth function $h$ vanishing on $\Sigma$ the product $h f $ is in $\cD^{\infty,\alpha}$ and $$  \Vert h f \Vert_{\cD^{k,\alpha}}\leq C_{k,h} \Vert \sigma \Vert_{k,\alpha}. $$
\item  $r \nabla r.\nabla f + \shalf f $ is in $\cD^{\infty,\alpha}$ and
$$  \Vert r \nabla r. \nabla f + \shalf f  \Vert_{\cD^{k,\alpha}}\leq C_{k} \Vert
\sigma \Vert_{k,\alpha}. $$
   \item With the choice of $W$ above:
 $$\tDelta f = \left( \nabla^{*}\nabla \sigma + \kappa \sigma\right)\zeta^{-\shalf} + \rho $$
where  $\rho\in \cD^{\infty,\alpha}$ and $$\Vert \rho\Vert_{\cD^{k,\alpha}}\leq C_{k} \Vert \sigma\Vert_{k+2,\alpha} . $$\end{enumerate}
\end{prop}

The first two items are straightforward and we only discuss the third item.
Here $\nabla^{*}\nabla$ is the usual covariant Laplacian on sections of $N^{1/2}$. The  essential point is that, expanding in the normal direction,  $\tDelta (\sigma \zeta^{-1/2})$ is $O(r^{-1/2})$, the leading term is given by the expression in the formula and all other terms are $O(r^{1/2})$. The proof involves a calculation in Fermi co-ordinates which we only give in outline. The standard formula for the Laplacian in local coordinates is
$$   \Delta_{g} f = g^{-\shalf} \left( g^{\shalf} g^{ij} f_{,j}\right)_{, i} $$
where $g={\rm det}(g_{ij})$. If we set $W= g^{-\squart}$ and $\tDelta = W^{-1}\Delta_{g} W$ then one finds that
\begin{equation}  \tDelta f = \left( g^{ij} f_{j}\right)_{,i} +  (W^{-1} \Delta W) f . \end{equation}
In our situation we have co-ordinates $x_{1}, x_{2}, t_{1}, \dots t_{n-2}$ where the $t_{i}$ are coordinates on $\Sigma$ and $x_{1}, x_{2}$ are in the normal direction. We choose the $t_{i}$ so that the volume form of $\Sigma$ is $dt_{1}\dots dt_{n-2}$, so the standard volume form on the total space of $N$ is $dx_{1}dx_{2}dt_{1}\dots dt_{n-2}$ and $W=V^{-1/2}=g^{-\squart}$. The function $f$ is ${\rm Re} ( \sigma(t) z^{-1/2})$ where $z=x_{i}+ i x_{2}$ and $\sigma$ is a complex-valued function of the $t_{j}$. One finds that, up to terms of order $r^{1/2}$,
$$  \left( g^{ij} f_{j}\right)_{,i} = (\nabla^{*}\nabla \sigma) z^{-1/2} - \squart K_{N} f . $$
For the other term in  (30) we need to compute the smooth function $\Delta W$. A calculation shows that on $\Sigma$ this is equal to $-\shalf {\rm Tr}_{N} {\rm Ric} + \sthreequart \vert \mu\vert^{2} $ and this gives the formula in the Proposition for the leading term.  (The  expansion of the volume form in Fermi coordinates is a well-studied topic, see for example \cite{kn:G}, Thm. 9.22.) The fact that the higher order term $\rho$ is in $\cD^{\infty,\alpha}$ and the bounds on the norms of $\rho$ are straightforward. 


\

\

To prove Theorem 4,  write $\tDelta (\sigma\zeta^{-1/2}) = \tau \zeta^{-1/2} + \rho$ where $\tau$ is given by the expression in item (3) of Proposition 10.  We consider $\tDelta ( r^{2} \tau \zeta^{-1/2})$ which can be written as
$$ \tDelta ( r^{2} \tau \zeta^{-\shalf}) =   (\Delta r^{2})\tau\zeta^{-\shalf} + 4 r \nabla r.\nabla(\tau \zeta^{-\shalf}) + r^{2}\tDelta (\tau \zeta^{-\shalf}) + 4\left( r W^{-1}\nabla W. \nabla r \right) \tau\zeta^{-\shalf}. $$
We have $\Delta r^{2} = 4$ on $\Sigma$ and the results of Proposition 10, applied to various sections of $N^{\shalf}$, imply that
$$    \tDelta ( \shalf r^{2} \tau \zeta^{-\shalf}) = \tau \zeta^{-\shalf} + \rho_{1}, $$
where
$$      \Vert \rho_{1} \Vert_{\cD^{k,\alpha}}\leq C_{k} \Vert \tau\Vert_{k+2,\alpha}. $$
 
 Now we can solve the equation $\tDelta u = \rho-\rho_{1}$ with $u\in \cC^{\infty,\alpha}$. So $\tDelta (\sigma \zeta^{-\shalf} - \shalf r^{2}\tau \zeta^{-\shalf} - u )=0$ and we have our solution
$$  \tQ= (\sigma \zeta^{-1/2} - \shalf r^{2}\tau \zeta^{-1/2} - u ) . $$

Uniqueness follows from the easy  fact that $r^{2}\tau \zeta^{-\shalf}$ is in $\cL^{2}_{1}$.

The solution $\tQ$ has an expansion
$$   \tQ= \sigma\zeta^{-\shalf} + \cA \zeta^{\shalf} + O(r^{\sthreetwo}) $$
where $\cA=A(u_{\sigma})\in \Gamma(\Sigma, N^{-\shalf})$. We define an operator
$$  P: \Gamma(\Sigma, N^{\shalf})\rightarrow \Gamma(\Sigma, N^{-\shalf})$$
by $P(\sigma)= \cA$. Combining (16) and the last estimate in Theorem 4, we see that $P$ is a bounded operator \begin{equation}P:C^{k+4,\alpha}\rightarrow C^{k+1,\alpha+\shalf}\end{equation} for each $k$. In fact, as we shall see later, this is not optimal and $P$ maps $C^{k+4,\alpha}$ to $C^{k+3,\alpha}$ but the difference is irrelevant for Nash-Moser theory.  

Notice that we get the same map $P$ for any choice of function $W$. We can express the same results in terms of the ordinary Laplace operator. 
\begin{cor}

For any $\sigma\in \Gamma(N^{-1/2})$ there is a unique section $Q$ of $E_{\bC}$ over $M\setminus \Sigma$ with $\Delta_{g} Q=0$ such that
$Q= \sigma\zeta^{-1/2}+ q$ where $q$ is in $\cL^{2}_{1}$. The section $Q$ satisfies 
$$  Q= \sigma\zeta^{-\shalf} + (P\sigma)\zeta^{\shalf} -\shalf {\rm Re}(\overline{\mu} \zeta) \sigma \zeta^{-\shalf} + O(r^{\sthreetwo}). $$ 
\end{cor}

The uniqueness is immediate from the Hilbert space theory and the existence follows from Theorem 4 by considering $Q=W \tQ$

\

The operator $P$ is analogous to the Dirichlet-to-Neumann operator for boundary value problems. Recall from the remark in Subsection 2.1 that in the model $\bC\times \bR^{n-2}$, if we restrict to functions of the form
$f=r^{-1/2} g(r,t)e^{i\theta/2}$, the Laplace equation  $\Delta f=0$ becomes  the ordinary Laplace equation for $g$ on a half-space in $\bR^{n+1}$. Globally, suppose that there is a circle action on $(M,g)$ with fixed point set $\Sigma$, so the quotient $M/S^{1}$ is an $(n-1)$-manifold $\Omega$ with boundary $\Sigma$. The harmonic sections of $E$ which transform with weight $\shalf$ under this action correspond to solutions of an equation of Laplace type on $\Omega$ and the operator $P$ is the usual Dirichlet-to-Neumann operator (for complex valued functions).

\subsection{The derivative formula}

We now compute the covariant derivative of our section $\bA$. Recall that is a linear map
$$  \frac{\delta \bA}{\delta\Sigma}:\Gamma(\Sigma, N)\rightarrow \Gamma(\Sigma, N^{-\shalf}). $$
We have natural pairings $N\otimes N^{-1/2}\rightarrow N^{1/2}$ and $N\otimes
N^{-3/2}\rightarrow N^{-1/2}$ and a section $\bB$ of $N^{-\sthreetwo}$.
\begin{prop} The derivative is given by
$$\left(\frac{\delta \bA}{\delta \Sigma}\right) v = \sthreetwo \bB v - \shalf P(\bA v) -\shalf \langle \mu,v\rangle \bA. $$ 
\end{prop}
 Here $\langle \mu,v\rangle$ is the {\it real} inner product of the vectors $\mu,v\in N$.

Assuming this, the proof of Theorem 3 is easy. If $\bB$ is nowhere-zero multiplication by $\stwothirds \bB^{-1}$ is an inverse
to the derivative up to $\bA$-quadratic errors 
$$  Q_{1}(\bA,v)=  
-\sthird \bB^{-1} P(\bA v)- \sthird \bB^{-1} \langle \mu,v\rangle \bA$$ and $$Q_{2}(\bA,\eta)=  -\sthird P(  \bA \bB^{-1} \eta)-\sthird \langle \mu,\bB^{-1} \eta\rangle A $$
 It is not hard to see that
these maps are smooth tame maps of the appropriate arguments.

\

We now prove Proposition 11. The statement is invariant under the diffeomorphisms of $M$, so it suffices to prove this with $(\Sigma, g)$ equal to our base point $(\Sigma_{0}, g_{0})$ in $\cS\times \cM$ over which we have a preferred point in $\cP$ given by the identity map. Starting with a section $v$ of $N$ we take the horizontal lift to a $1$-jet $\tilde{v}$ with $\nabla_{N} \tilde{v}=0$. We extend this $1$-jet to a vector field on $M$, which we also denote by $\tilde{v}$. Let $\lambda_{t}:M\rightarrow M$ be a $1$-parameter family of diffeomorphisms with derivative $\tilde{v}$ at $t=0$  and let $g_{t}= \lambda^{*}_{t} g$. We can suppose that $\lambda_{t}$ is chosen so that the metrics $g_{t}$ are compatible with the fixed normal structure on $\Sigma$. For each $t$ we have a solution $\phi_{t}$ of $\Delta_{g} \phi_{t}=0$ with singular set $\Sigma_{t}= \lambda_{t}(\Sigma)$. So the pull-back $\Phi_{t}=\lambda_{t}^{*}(\phi_{t})$  has singular set the fixed submanifold $\Sigma$ and is characterised as the solution of the equation $\Delta_{t} \Phi_{t}=0$. where $\Delta_{t}$ is the Laplacian of the metric $g_{t}$. By definition, what we want to compute is the time derivative 
\begin{equation}  \frac{\delta \bA}{\delta \Sigma}(v)=  \ddt A (\Phi_{t})= A(  \Phi'). \end{equation}
Where we have written $\Phi'=\ddt \Phi$.
(Here, and below, all $t$-derivatives are evaluated at $t=0$. )

{\bf Remark} Strictly speaking, we defined $A$ on the $\cC^{\infty,\alpha}$ sections whereas $\Phi_{t}$ is in $W_{t} \cC^{\infty,\alpha}$ for a $1$-parameter family of functions $W_{t}$ depending on the mean curvature of $\Sigma$ in the metrics $g_{t}$. Thus we should write
$$  \ddt (A(W_{t}^{-1} \Phi_{t})). $$
However $W_{t}$ does not affect the leading term in the asymptotic expansion and we can safely ignore it here.

 Let $\Delta'$ be the derivative of the family of operators $\Delta_{t}$ at $t=0$.  If $u$ is any section then
$$   \Delta_{\lambda^{*}_{t} g} \lambda^{*}_{t} u = \lambda_{t}^{*} (\Delta_{g} u). $$
Differentiating this, we get
$$  \Delta'(u)= \nabla_{\tilde{v}}\Delta u - \Delta\left( \nabla_{\tilde{v}} u\right). $$

On the other hand, differentiating the equation
$\Delta_{t} \Phi_{t}=0$ gives
$$  \Delta_{g}( \Phi')= -\Delta' (\phi). $$
So we have

$$\Delta_{g}( \Phi' -\nabla_{v} \phi) =0 $$ on $M\setminus \Sigma$. 

Now $ W_{t}\Phi_{t}$ is a smooth path in $\cC^{\infty,\alpha}$. This implies that $\Phi_{t}$ is a smooth path in $\cD^{\infty,\alpha}$ and hence $\Phi'$ is in $\cD^{\infty,\alpha}$ and thus  also in $ \cL^{2}_{1}$. We know that $\phi$ has an asymptotic description
$$  \phi= {\rm Re}\left( A \zeta^{\shalf} + B \zeta^{\sthreetwo}\right)  -\shalf {\rm Re}(\overline{\mu}\zeta) {\rm Re} (A \zeta^{1/2})+  O(r^{\sfivetwo}), $$
and the derivative of the $O(r^{\sfivetwo})$ term is $O(r^{\sthreetwo})$. It is straightforward to see that the condition $\nabla_{N}\tilde{v}=0$ implies that
\begin{equation}  \nabla_{\tilde v} \left( A \zeta^{\shalf}\right)=  \shalf (Av) \zeta^{-\shalf}) +O(r^{\sthreetwo}), \end{equation}
that
\begin{equation}
\nabla_{\tilde v} \left( B \zeta^{\sthreetwo}\right)=  \sthreetwo (Bv) \zeta^{\shalf} + O(r^{\sfivetwo})\end{equation}
and we have 
$$
\nabla_{\tilde{v}}\  {\rm Re}( { \overline{\mu} \zeta}) = \langle v,\mu\rangle $$
on $\Sigma$. From this we see that $ \nabla_{\tilde{v}} \phi$ differs from $\shalf {\rm Re} (Av) \zeta^{-\shalf}$ by an element of $\cL^{2}_{1}$, so by Corollary 1 we have $ \Phi'-\nabla_{\tilde v}\phi= -\shalf{\rm Re}\ Q(Av)$. That is
\begin{equation}
\Phi'-\nabla_{\tilde{v}}\phi= -\shalf {\rm Re}\left( Av \zeta^{-\shalf} + P(Av) \zeta^{\shalf}- \shalf {\rm Re}(\overline{\mu}\zeta) Av \zeta^{-\shalf}\right) +O(r^{\sthreetwo}).  \end{equation}
Comparing the $O(r^{\shalf})$ terms, we see that
$$  \Phi'= {\rm Re}\left( \left(\sthreetwo Bv -\shalf P(Av) - \shalf \langle \mu,v\rangle A\right) \zeta^{\shalf}\right) + O(r^{3/2}),$$
which gives the formula stated for the leading term $A(\Phi')$.

\subsection{Further analysis of the operator $P$}

We begin with the Greens function $G$ of the Laplacian $\Delta$ on $\Gamma(E)$. For distinct  points $p,q\in M$ we have $G(p,q)\in E_{p}\otimes E_{q}$. We use the metric on $E$ to identify
  $$  E_{p}\otimes E_{q}= {\rm Hom}(E_{p}, E_{q})= {\rm Hom}(E_{q}, E_{p}), $$
and $G(p,q)$ is symmetric in $p,q$. Now fix $q\in M\setminus \Sigma$ so
$G(\ ,q)$ is harmonic away from $q$ and has a leading term of order $r^{1/2}$ near $\Sigma$. For $t_{1}\in \Sigma$ this leading term gives $H(t_{1},q)\in N_{t}^{-1/2}\otimes_{\bR} E_{q}$, so for any section $u$ with $\Delta u$ in $C^{,\alpha}$ the leading term $A=A(u)\in \Gamma(N^{-1/2})$ is 
$$  A(u)(t_{1})= \int_{M} H(t_{1}, q) (\Delta u)(q) dq. $$
For fixed $t_{1}$ the section $H(t_{1},\ )$ is harmonic over $M\setminus \{t_{1}\}$ so for $t_{2}\neq t_{1}$ there is a leading asymptotic term
$$\Gamma(t_{1}, t_{2})\in N_{t_{1}}^{1/2}\otimes_{\bR} N_{t_{2}}^{-1/2}. $$

For $t_{2} $ close to $t_{1}$ let $\bone$ be the element of $N_{t_{1}}^{1/2}\otimes N_{t_{2}}^{-1/2}$ defined by parallel transport of the identity along the minimal geodesic from $t_{1}$ to $t_{2}$. As convenient notation, we write $\vert t_{1}-t_{2}\vert$ for the Riemannian distance between $t_{1}, t_{2}$. 
\begin{thm}\begin{enumerate}
\item For $t_{1}$ close to $t_{2}$ in $\Sigma$ the section $\Gamma$ has asymptotic behaviour 
$$  \Gamma(t_{1}, t_{2}) = \kappa_{n}\vert t_{1}- t_{2}\vert^{1-n}\bone  + O(\log \vert t_{1}-t_{2}\vert \vert t_{1}-t_{2}\vert ^{3-n}), $$ where $\kappa_{n}$ is the constant defined in Subsection 2.1.
\item
For $\sigma\in \Gamma(N^{1/2})$
$$  P\sigma(t_{1}) = \pi \lim_{\delta\rightarrow 0}\left( \int_{\vert t_{1}-t_{2}\vert \geq \delta} \Gamma (t_{1}, t_{2}) \sigma(t_{2}) dt_{2} - {\rm Vol}(S^{n-3}) \kappa_{n} \delta^{-1} \ \sigma(t_{1}) \right) + (\overline{\mu}\sigma)(t_{1}).$$

\end{enumerate}
\end{thm}

{\bf Remark} 
The formula in Theorem 5  and general theory show that the operator $P$ is a pseudodifferential operator of order $-1$, with the same symbol as the square root of the Laplace operator on $\Gamma (\Sigma, N^{-1/2})$. This is not surprising in view of the connection with the Dirichlet-to-Neumann operators mentioned at the beginning of this section. This suggests that it is not possible to prove our main result using the more elementary implicit function theorem in Banach spaces. We can choose function spaces such that
$\frac{\partial \bA}{\partial \Sigma}$ is bounded when $\bA\neq 0$: for example we can consider it as a map from $C^{k,\alpha}$ to $C^{k-1,\alpha}$. But with such a choice of function spaces the derivative is not invertible when $\bA=0$. 

\

\

Fix $\sigma$ and let $Q$ be the solution of $\Delta Q=0$ with leading term
$\sigma\zeta^{-1/2}$ as in Corollary 1. Fix $q\in M\setminus \Sigma$ and for small $\epsilon>0$ let $M_{\epsilon}$ be the complement of the $\epsilon$-neighbourhood of $\Sigma$. We can suppose that $q$ is in $M_{\epsilon}$ and apply Green's formula on the manifold-with-boundary $M_{\epsilon}$:
\begin{equation}  Q(q)= \int_{\partial M_{\epsilon}}  \frac{\partial G}{\partial \nu}Q - \frac{\partial Q}{\partial \nu} G. \end{equation}

Here $G=G(\ ,q)$ and $\frac{\partial}{\partial \nu}$ denotes the normal derivative on $\partial M_{\epsilon}$. Using the asymptotic descriptions of $Q$ and $G$, a simple calculation identifies the limit of the right hand side of (36) as $\epsilon$ tends to $0$,  and we get
\begin{equation} Q(q)= \pi \int_{\Sigma} H(t_{2},q) \sigma(t_{2}) dt_{2}. \end{equation} 

Now if $t_{1}$ is outside the support of $\sigma$ it follows immediately from the definitions that:
$$  P(\sigma)(t_{1}) = \pi \int_{\Sigma} \Gamma(t_{1}, t_{2})\sigma(t_{2}) dt_{2}, $$
which is a special case of the second item in Theorem 5.

 To establish the second item of Theorem 5 in general we need to study the behaviour of $\Gamma$ near the diagonal and  particular the asymptotic behaviour stated in the first item of Theorem 5. The main point here is that there is no $O(\vert t_{1}-t_{2}\vert^{2-n})$ term in the expansion.

To simplify notation we will restrict attention to the case when  $n=3$, so $\Sigma$ is a $1$-dimensional manifold. The higher dimensional case is essentially the same.   Fix a point $\tau$ in $\Sigma$ and take standard coordinates $(z,t)$ centred at this point. For $p=(z,t)$ write $H(\tau,p) = h(z,t)$ and for $\tau'= (0,t)$ in $\Sigma$ write $\Gamma(\tau, \tau')= \gamma(t)$. Thus the second item of Theorem 5 is the assertion that 
\begin{equation}   \gamma(t) = \pi^{-1} t^{-2} + O\left(\log \vert t \vert)\right) \end{equation}
as $t\rightarrow 0$. 

Let $h_{0}(z,t)$ be the section in the flat model denoted by $H(z,t)$ in subsection 2.1 above. This has homogeneity $-3/2$ and, away from the origin,  satisfies the equation $\Delta_{0} h_{0}=0$ for the Euclidean Laplacian $\Delta_{0}$, whereas $h(z,t)$ satisfies the equation $\Delta_{g} h=0 $. We want to obtain an asymptotic expansion of $h$, starting with the flat model $h_{0}$, following standard procedures. This will have the form  
\begin{equation}  h(z,t)\sim h_{0}(z,t)+ h_{1}(z,t)+ \dots.  \end{equation}
 Write 
$$ \Delta_{g}= \Delta_{0}+ L_{1} + L_{2} + \dots, $$
where $L_{a}$ is a differential operator containing terms $p D$ where
$D$ is a constant co-efficient operator of order $d_{1}$ and $p$ is a polynomial of degree $d_{2}$ and $d_{2}-d_{1}=a$. So $L_{a}$ maps a function of homogeneity $\lambda$ to  one of homogeneity $\lambda-a$ while $\Delta_{0}$ maps to homogeneity $\lambda-2$. We construct the expansion by solving  equations for $\Delta_{0} h_{k}$ to successively remove the error terms. 

In spherical polar coordinates, with $R=\vert(z,t)\vert$, consider a section $R^{\lambda} f$ where $f$ is a section of a flat bundle over $S^{2}$ minus the poles.
Then $$  \Delta_{0}(R^{\lambda} f)= R^{\lambda-2}(\Delta_{S^{2}}+ \lambda(\lambda+1)) f).$$ 
So if $-\lambda(\lambda+1)$ is not an eigenvalue of $\Delta_{S^{2}}$ (acting on sections of this flat bundle) the operator $\Delta_{0}$ maps sections of homogeneity $\lambda$ onto those of homogeneity $\lambda-2$. If we do not encounter any such eigenvalues we can construct the asymptotic series using terms of the form $R^{\lambda} f_{\lambda}$.  If we encounter eigenvalues we may need to include terms of the form $(\log R)^{\mu} R^{\lambda} f_{\lambda,\mu} $. In any event it is standard to show that there is an asymptotic solution. Taking sufficiently many terms we get a finite sum
$$  \tilde{H}(z,t)=h_{0}+ \dots h_{k}, $$
with $\Delta_{g}\tH = O(R^{1/2})$.  Then $\Delta_{g} \tH\in C^{,\alpha}$ and we can use our global theory to see that $h- \tilde{H}$ is in $W \cC^{2,\alpha}$.
For each $i$ the section $h_{i}(z,t)$ has $h_{i}(z,t)\sim {\rm Re} (\gamma_{i}(t) z^{1/2})$ as $z\rightarrow 0$. If we write $\tilde{\gamma}(z,t)= \gamma_{0}+ \dots\gamma_{k}$ our theory shows that  
$\gamma-\tilde{\gamma}$ is in $C^{1,\alpha+1/2}$.

Given the existence of this asymptotic expansion, for the purposes of the proof of Theorem 5 we only need to study the first two terms $h_{0},h_{1}$ in the series. The leading term $h_{0}$ is of homogeneity $-3/2$ and we have an explicit formula (5).  When $\lambda=-1/2$ the expression $-\lambda(\lambda+1)$ is positive and hence not an eigenvalue,  thus we know that we can find a second term $h_{1}$ of homogeneity $-1/2$ to remove the $O(R^{-5/2})$ term in $\Delta_{g} h_{0}$.  

The explicit formula (5) shows that
$$  \gamma_{0}(t) =  \kappa_{3} t^{-2} = \pi^{-1} t^{-2} $$
 Thus for fixed $t\neq 0$ we have
$$  h_{0}(z,t)\sim \pi^{-1} z^{1/2} t^{-2} $$ as $z\rightarrow 0$. More precisely, using the scaling behaviour,  we have that for all $\vert z\vert\leq t$
\begin{equation}  \vert h_{0}(z,t)-\pi^{-1} z^{1/2} t^{-2}\vert\leq C \vert z \vert^{3/2} \vert t\vert^{-3}. 
\end{equation} 

Our next step is to identify $h_{1}$ explicitly in terms of $h_{0}$. 
The only terms in the differential operator $L_{1}$ come from the mean curvature of $\Sigma$. To simplify notation we can suppose our coordinates are chosen with $z=x_{1}+ i x_{2}$ and the mean curvature at the point $\tau$ equal to
$$  \mu(\tau) = m \frac{\partial}{\partial x_{1}}.$$
Then 
$$  L_{1}f =  m\left( \frac{\partial f}{\partial x_{1}} - 2 \frac{\partial^{2} f}{\partial t^{2}}\right) $$

The first term is the term we encountered in Proposition 6.   To find $h_{1}$ we need to solve the equation $\Delta_{0} h_{1}= - L_{1} (h_{0})$.   Using $\Delta_{0} h_{0}=0$ we get
$$  \Delta_{0} \left( \frac{x_{1} h_{0}}{2}\right)= \frac{\partial h_{0}}{\partial x_{1}} $$ and 
$$  \Delta_{0}(\frac{x_{1}^{2}-x^{2}_{2}}{2} \frac{\partial h_{0}}{\partial x_{1}}+ x_{1} x_{2} \frac{\partial h_{0}}{\partial x_{2}}) = -2x_{1} \frac{\partial^{2} h_{0}}{\partial t^{2}}. $$
It follows that 
\begin{equation}
h_{1}= m \left(\frac{x_{1}^{2}-x^{2}_{2}}{2} \frac{\partial h_{0}}{\partial
x_{1}}+ x_{1} x_{2} \frac{\partial h_{0}}{\partial x_{2}} + x_{1} h_{0}\right) \end{equation}

The $O(z^{1/2})$ term in $h_{1}$ clearly vanishes so $\gamma_{1}=0$. Thus we have the first item in Theorem 5:
\begin{equation} \gamma(t)= \pi^{-1} t^{-2} + O(\log \vert t\vert). \end{equation}
  Note that for the next term in the series, with $\lambda=1/2$, we do encounter an eigenvalue (since $z^{\shalf}$ is harmonic) and we expect to introduce   a term $h_{2}= O(R^{\shalf} \log R)$ and an   $O(\log \vert t\vert ) $ term in $\gamma$.

For fixed $z$ the explicit formula (5) shows that the integral
$$  I_{0}(z)= \int_{-\infty}^{\infty} h_{0}(z,t) dt, $$
is defined and one calculates
$$  I_{0}(z)=  z^{-1/2}. $$
 In a similar way,  we see that the integral
$$  I_{1}(z)= \int_{-\infty}^{\infty} h_{1}(z,t) dt, $$
is defined and  our formula (41) for $h_{1}$ implies that
$$ I_{1}= m \left( \frac{x_{1}^{2}-x^{2}_{2}}{2} \frac{\partial I_{0}}{\partial x_{1}} + x_{1}x_{2}\frac{\partial I_{0}}{\partial x_{2}} + x_{1} I_{0}\right).
$$ We conclude that
\begin{equation} I_{1}(z)= \frac{-m}{4} z^{1/2} + m x_{1} z^{-1/2}.\end{equation}

With these preparations in place we can move on to prove the second item in Theorem 5. Writing $\sigma_{0}$ for $\sigma(\tau)$ and $P_{0}$ for
$(P\sigma)(\tau)$, we know that
\begin{equation} Q(z,0)= \sigma_{0} z^{-\shalf} -\shalf m x_{1} \sigma_{0}z^{-\shalf}+  P_{0} z^{\shalf} +O(r^{\sthreetwo}). \end{equation}
On the other hand we have from (37):
\begin{equation} Q(z,0)= \pi \int_{\Sigma} H(t, (z,0))\sigma(t) dt. \end{equation}
Let $\delta$ be a small number and $r=\vert z\vert\leq \delta$. We write the integral in (45)  as $J_{1}+J_{2}$ with
$$  J_{1} =\pi \int_{\vert t\vert>\delta}  H(t, (z,0))\sigma(t) dt,
$$
and
$$ J_{2} = \pi \int_{\vert t\vert<\delta}  H(t, (z,0))\sigma(t)
dt. $$
(Here we understand that $J_{1}$ contains also the contribution away from our co-ordinate neighbourhood.)
We will take $z\rightarrow 0$, with $\delta$ fixed. Clearly as $z\rightarrow 0$
\begin{equation}    J_{1} \sim \pi \int_{\vert t\vert \geq \delta} (\gamma(t) \sigma(t) z^{1/2} )dt, \end{equation}
which is one of the terms appearing in the formula of Theorem 5. To analyse $J_{2}$  we use the asymptotic description of $H$ discussed above. Two observations simplify the calculations. First if we replace $\sigma(t)$ by the fixed value $\sigma_{0}$ we introduce an error of order $\delta r^{1/2}$. Second, in our asymptotic analysis above we considered $h(z,t)= H(0,(z,t)$ whereas the quantity appearing in (45) is $H(t,(z,0))$. It is straightforward to see that, up to an error of order $\delta r^{1/2}+ r^{3/2}$, we can interchange the two points of view and approximate $H(t,(z,0))$ by $h(z,-t)$. The conclusion is that
\begin{equation} J_{2}= \tilde{J}_{2} + O(\delta r^{1/2} + r^{3/2}) \end{equation}
where 
\begin{equation}  \tilde{J}_{2} = \pi \int_{-\delta}^{\delta}  h(t,z) \sigma_{0} dt. \end{equation}
Now write
$$ \tilde{J}_{2} = \tilde{J}_{3}+\tilde{J}_{4}$$
where
$$ \tilde{J}_{3}= \pi \int_{-\delta}^{\delta} (h_{0} + h_{1})(z,t)\sigma_{0} dt, $$
and
$$  \tilde{J}_{4}= \pi \int_{-\delta}^{\delta} ( h- (h_{0}+ h_{1}))(z,t)\sigma_{0} dt. $$

The most important term is $\tilde{J}_{3}$. For $a=0$ or $1$, write
$$  I_{a,\delta}(z) = \int_{-\delta}^{\delta} h_{a}(z,t). $$
Then (40) implies that
$$  I_{0,\delta}(z) = I_{0}(z) - 2 \delta^{-1} \pi^{-1} z^{1/2} + O(\delta^{-2} r^{3/2}).$$
Similarly
  $$  I_{1,\delta}(z)= I_{1}(z) + O(\delta^{-1} r^{3/2}). $$
  So we get
  $$  \tilde{J}_{3}= \left(\pi I_{0}(z)+I_{1}(z)-2\delta^{-1}  z^{1/2}\right) \sigma_{0}  +O(\delta^{-2}r^{3/2}). $$
The main contribution to $\tilde{J}_{4}$ comes from the next term $h_{2}$ in the asymptotic expansion and one finds that
$$  \tilde{J}_{4}= O(r^{1/2} \delta \log \delta). $$
Putting this together we get
\begin{equation}  J_{2}= \left(\pi I_{0}(z)+ \pi I_{1}(z) -2\delta^{-1}\ z^{\shalf}\right)\sigma_{0} + O(r^{\shalf}\delta \log \delta )+ O(\delta^{-2} r^{\sthreetwo}). \end{equation}
We now compare $J_{1}+J_{2}$ with (44) and take $z\rightarrow 0$. We use (49) and the formulae above for $I_{0}, I_{1}$. The $O(r^{\sthreetwo})$ terms on either side match up (as  they must). The $O(r^{\shalf})$ terms of the form $m x_{1} z^{-1/2}$ also match up. Comparing the $z^{\shalf}$ terms we get
$$   P_{0}= \pi \int_{\vert t\vert\geq \delta} \gamma(t)\sigma(t) dt - 2 \delta^{-1} \sigma_{0} + m\sigma_{0} +O(\delta \log \delta), $$
which gives the formula in the second item of Theorem 5.



\end{document}